%% file: 2001-13.tex
\documentclass{gtart}

\input gtoutput

\lognumber{135}
\volumenumber{5}\papernumber{13}\volumeyear{2001}
\pagenumbers{369}{398}
\accepted{8 April 2001}
\received{12 September 2000}
\published{20 April 2001}
\proposed{Walter Neumann}
\seconded{Cameron Gordon, David Gabai}

\usepackage{rlepsf}
\let\relabela\adjustrelabel

\usepackage{amsmath}
\usepackage{amsfonts}
\usepackage{amssymb}
\usepackage{latexsym}

\newcommand{\begriff}[1]{{\em #1}}

\newcommand{\x}{\mathfrak x}
\renewcommand{\d}{\partial}
\newcommand{\nach}{\to}
\renewcommand{\t}{\mathcal T}

\newtheorem{theorem}{Theorem}
\newtheorem{proposition}{Proposition}
\newtheorem{lemma}{Lemma}
\newtheorem{corollary}{Corollary}

\theoremstyle{definition}
\newtheorem{definition}{Definition}
\newtheorem{construction}{Construction}

\theoremstyle{remark}
\newtheorem{remark}{Remark}

\begin{document}
\title{The size of triangulations supporting a given link}  
\author{Simon A King}
\address{Institut de Recherche Math\'ematique Avanc\'ee\\Strasbourg,
  France}
\asciiaddress{Institut de Recherche Mathematique Avancee, Strasbourg,
  France}
\email{king@math.u-strasbg.fr}
\begin{abstract}
  Let $\t$ be a triangulation of $S^3$ containing a link $L$ in its
  1--skeleton. We give an explicit lower bound for the number of
  tetrahedra of $\t$ in terms of the bridge number of $L$. Our proof
  is based on the theory of almost normal surfaces.
\end{abstract}
\asciiabstract{
  Let T be a triangulation of S^3 containing a link L in its
  1-skeleton. We give an explicit lower bound for the number of
  tetrahedra of T in terms of the bridge number of L. Our proof
  is based on the theory of almost normal surfaces.}


\asciikeywords{Link, triangulation, bridge number, Rubinstein-Thompson
algorithm, normal surfaces}

\maketitlepage

\section{Introduction}
\label{sec:intro}

In this paper, we prove the following result.
\begin{theorem}\label{thm:main}
  Let $L\subset S^3$ be a tame link with bridge number $b(L)$.
  Let $\t$ be a triangulation of $S^3$ with $n$ tetrahedra such
  that $L$ is contained in the 1--skeleton of $\t$. 
  Then $$ n >\frac{1}{14} \sqrt{\log_2 b(L)},$$ 
  or equivalently $$b(L) < 2^{196 n^2}.$$
\end{theorem}
The definition of the bridge number
can be found, for instance, in~\cite{burde}. So far as is known to the
author, Theorem~\ref{thm:main} gives the first  estimate for $n$ in
terms of $L$ that does not rely on additional geometric or combinatorial
assumptions on $\t$.  
We show in~\cite{king} that the bound for $b(L)$ in
Theorem~\ref{thm:main} can not be replaced by a sub-exponential bound in
$n$. 
More precisely, there is a constant $c\in\mathbb R$ such that for any
$i\in \mathbb N$ 
there is a triangulation $\t_i$ of $S^3$ with $\le c\cdot i$
tetrahedra, containing a two-component link $L_i$ in its 1--skeleton
with $b(L_i)> 2^{i-2}$. 

The relationship of geometric and combinatorial properties of a
triangulation of $S^3$ with the knots in its 1--skeleton 
has been studied earlier, see~\cite{goodrick},
\cite{lickorish},
\cite{armentrout}, \cite{ehrenborghachimori},~\cite{hachimoriziegler}.
For any knot $K\subset S^3$ there is a triangulation of $S^3$ such
that $K$ is formed by three edges, see~\cite{furch}. 
Let $\t$ be a triangulation of $S^3$ with $n$ tetrahedra and let 
$K\subset S^3$ be a knot formed by a path of $k$ edges. If $\t$ is
shellable (see~\cite{ehrenborghachimori}) or
the dual cellular decomposition is shellable (see~\cite{armentrout}), then
$b(K)\le \frac 12 k$. If $\t$ is vertex decomposable then $b(K)\le \frac
13 k$, 
see~\cite{ehrenborghachimori}.

We reduce Theorem~\ref{thm:main} to Theorem~\ref{prop:} below, for
which we need some definitions. Denote $I=[0,1]$. Let $M$ be a
closed 3--manifold with a triangulation $\t$. The $i$--skeleton of $\t$ is
denoted by $\t^i$. 
Let $S$ be a surface and let $H\co S\times I\to M$ be an embedding, so that
$\t^1\subset H(S^2\times I)$.
A point $x\in \t^1$
is a \emph{critical point} of $H$ if 
$H_\xi = H(S\times \xi)$ is not transversal to $\t^1$ in $x$, for some $\xi \in I$. 
We call $H$ a \emph{$\t^1$--Morse embedding}, if $H$ is in general
position with 
respect to $\t^1$;
we give a more precise definition in Section~\ref{sec:almost}. 
Denote by $c(H,\t^1)$ the number of critical points of $H$. 

\begin{theorem}\label{prop:}
  Let $\t$ be a triangulation of $S^3$ with $n$ tetrahedra.
  There is a $\t^1$--Morse embedding $H\co S^2\times I\nach S^3$ such that 
  $\t^1\subset H(S^2\times I)$ and $c(H,\t^1)<2^{196 n^2}$.
\end{theorem}
For a link $L\subset \t^1$, it is easy to see that $b(L)\le\frac 12
\min_H \{c(H,\t^1)\}$, where the minimum is taken over all $\t^1$--Morse
embeddings $H\co S^2\times I\to S^3$ with $L\subset H(S^2\times I)$. Thus
Theorem~\ref{thm:main} is a corollary of Theorem~\ref{prop:}. 

Our proof of Theorem~\ref{prop:} is based on the theory of almost
2--normal surfaces. Kneser~\cite{kneser} introduced 1--normal surfaces
in his study of connected sums of 3--manifolds. 
The theory of 1--normal surfaces was further developed by Haken
(see~\cite{haken1}, \cite{haken2}). It led to
a classification algorithm for knots and for sufficiently
large 3--manifolds, see for instance~\cite{hemion},~\cite{matveev2}. 
The more general notion of almost 2--normal surfaces is due
to Rubinstein~\cite{rubinstein}. With this 
concept, Rubinstein and Thompson found a recognition algorithm for
$S^3$, see~\cite{rubinstein}, \cite{thompson},~\cite{matv}. 
Based on the results discussed in a preliminary version of this
paper~\cite{king2}, the author~\cite{king} and 
Mijatovi\'c~\cite{mijatovic} independently obtained a recognition
algorithm for $S^3$ using local transformations of triangulations.

We outline here the proof of Theorem~\ref{prop:}. Let $\t$ be a
triangulation of $S^3$ with $n$ tetrahedra.  
If  $S\subset S^3$ is an embedded surface and $S\cap \t^1$ is finite,
then set $\|S\|= \text{card}(S\cap \t^1)$. Let 
$S_1,\dots, S_k\subset S^3$ be surfaces. A surface that is 
obtained by joining $S_1,\dots, S_k$ with some small 
tubes in $M\setminus \t^1$ is called a \emph{tube sum} of $S_1,\dots,
S_k$. 

Based on the Rubinstein--Thompson algorithm, we construct a  system
$\tilde\Sigma\subset S^3$ of pairwise disjoint 2--normal 2--spheres such
that  $\|\tilde\Sigma\|$ is bounded in terms of $n$ and any
1--normal sphere in $S^3\setminus \tilde\Sigma$ is parallel to a connected
component of $\tilde\Sigma$.  The bound for $\|\tilde\Sigma\|$ can
be seen as part of a complexity analysis for the Rubinstein--Thompson
algorithm and relies on results on integer programming.

A $\t^1$--Morse embedding $H$ then is constructed ``piecewise'' in the
connected 
components of $S^3\setminus \tilde\Sigma$, which means the
following. There are numbers $0<\xi_1<\dots < \xi_m< 1$ such that: 
\begin{enumerate}
\item $\|H_0\| = \|H_1\|= 0$.
\item There is one critical value of $H|[0,\xi_1]$, corresponding
  to a vertex $x_0\in\t^0$. The set of critical points
  of $H|[\xi_m,1]$ is $\t^0\setminus\{x_0\}$.   
\item For any $i=1,\dots,m$, the sphere $H_{\xi_i}$ is a tube sum of
  components of $\tilde\Sigma$.
\item The critical points of $H|[\xi_i,\xi_{i+1}]$ are contained in a
  single connected component $N_i$ of $S^3\setminus \tilde\Sigma$.
\item The function $\xi\mapsto\|H_\xi\|$ is monotone in any interval
  $[\xi_i,\xi_{i+1}]$, for $i=1,\dots, m-1$. 
\end{enumerate}
\begin{figure}[ht!]
\cl{\relabelbox\small
\epsfbox{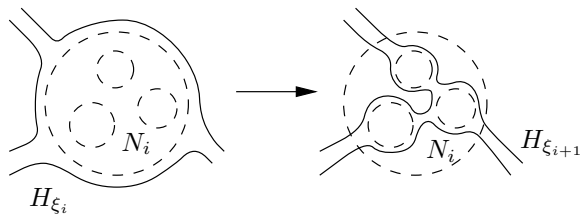}
\relabel {N}{$N_i$}
\relabel {N1}{$N_i$}
\relabel {H1}{$H_{\xi_i}$}
\relabel {H}{$H_{\xi_{i+1}}$}
\endrelabelbox}
\caption{About the construction of $H$}
    \label{fig:thatsH}
\end{figure}
This is depicted in Figure~\ref{fig:thatsH}, where the components of
$\tilde \Sigma$ are dotted. The components $N_i$ run over all
components of $S^3\setminus \tilde\Sigma$ that are not regular
neighbourhoods of vertices of $\t$. Thus an estimate for $m$ is
obtained by an estimate for the number of components of $\tilde
\Sigma$. By monotonicity of $\|H_\xi\|$, the number of critical points
in $N_i$ is bounded by $\frac 12\, \|\d N_i\|\le \frac 12\,
\|\tilde\Sigma\|$. This together with the bound for $m$ yields the
claimed estimate for $c(H,\t^1)$.

The main difficulty in constructing $H$ is to assure property (5). For
this, we introduce the notions of upper and lower reductions. If $S'$
is an upper (resp.\ lower) reduction of a surfaces $S\subset S^3$,
then $S$ is isotopic to $S'$ such that $\|\cdot\|$ is monotonely
non-increasing under the isotopy.
Let $N$ be a connected component of $S^3\setminus \tilde\Sigma$ with
$\d N=S_0\cup S_1\cup\dots \cup S_k$. 
We show that there is a tube sum $S$ of $S_1,\dots,S_k$ such that
either $S$ is a lower reduction of $S_0$, or $S_0$ is an upper reduction
of $S$. 
Finally, if $H_{\xi_i}$ is a tube sum of $S_0$ with some surface
$S'\subset S^3\setminus N$, then $H|[\xi_i,\xi_{i+1}]$ is induced by
the lower reductions (resp.\ the inverse of the upper reductions)
relating $S_0$ with $S$. Then 
$H_{\xi_{i+1}}$ is a tube sum of $S$ with $S'$, assuring properties
(3)--(5).

The paper is organized as follows. 
In Section~\ref{sec:knormal}, we recall basic properties of 
$k$--normal surfaces. 
It is well known that the set of 1--normal surfaces in a triangulated
3--manifold is additively generated by so-called \emph{fundamental
  surfaces}. In Section~\ref{sec:fundamental}, we generalize this to
2--normal surfaces contained in \emph{sub-manifolds} of  triangulated
3--manifolds.
The system $\tilde\Sigma$ of 2--normal spheres is constructed in
Section~\ref{sec:constrN}, in the more
general setting of closed orientable 3--manifolds.
In Section~\ref{sec:almost}, we recall the notions of almost
$k$--normal surfaces (see~\cite{matv}) and of impermeable surfaces
(see~\cite{thompson}), and introduce the new  notion of
split equivalence. We discuss the close relationship of almost
2--normal surfaces and impermeable surfaces. This relationship is well
known (see~\cite{thompson},~\cite{matv}), but the proofs are only
partly available. For completeness we give a proof in the last
Section~\ref{sec:imperm-a2}.
In Section~\ref{sec:1-2-normal} we exhibit some useful properties of
almost 1--normal surfaces. 
The notions of upper and lower reductions are introduced in
Section~\ref{sec:reduction}.  
The proof of
Theorem~\ref{prop:} is finished in Section~\ref{sec:theproof}. 

In this paper, we denote by $\#(X)$ the number of connected components
of a topological space $X$. 
If $X$ is a tame subset of a 3--manifold $M$, then $U(X)\subset M$
denotes a regular neighbourhood of $X$ in $M$. For a
triangulation $\t$ of $M$, the number of its tetrahedra is denoted by
$t(\t)$. 

\rk{Acknowledgements}
I would like to thank Professor Sergei V Matveev and my scientific supervisor
Professor Vladimir G Turaev for many interesting discussions and for
helpful comments on this paper.


\section{A survey of $k$--normal surfaces}
\label{sec:knormal}

Let $M$ be a closed 3--manifold with a
triangulation $\t$. The number of its tetrahedra is denoted by $t(\t)$. 
An \begriff{isotopy mod $\t^n$} is an ambient isotopy of $M$ that
fixes any simplex of $\t^n$ set-wise. Some authors call an isotopy mod
$\t^2$ a normal isotopy.

\begin{definition}
  Let $\sigma$ be a 2--simplex and let $\gamma\subset \sigma$ be a
  closed embedded 
  arc with $\gamma\cap\d\sigma=\d\gamma$, disjoint to the vertices of
  $\sigma$. If $\gamma$ connects two 
  different edges of $\sigma$ then $\gamma$ is called a
  \begriff{normal arc}. Otherwise, $\gamma$ is called a
  \begriff{return}. 
\end{definition}
We denote the number of connected components of a topological space
$X$ by $\#(X)$.
Let $\sigma$ be a 2--simplex with edges $e_1,e_2,e_3$. If
$\Gamma\subset\sigma$ is a system of 
normal arcs, then $\Gamma$ is determined by $\Gamma\cap\d
\sigma$, up to isotopy constant on $\d \sigma$, and
$e_1$ is connected with $e_2$ by
$\frac 12\left(\#(\Gamma\cap e_1) + \#(\Gamma\cap e_2) - \#(\Gamma\cap
e_3)\right)$ arcs in $\Gamma$. 
\begin{definition}
  Let $S\subset M$ be a closed embedded surface transversal to $\t^2$. 
  We call $S$ \begriff{pre-normal}, if $S\setminus\t^2$ is a disjoint
  union of discs and $S\cap\t^2$ is a union of normal arcs in the
  2--simplices of $\t$. 
\end{definition}
The set $S\cap \t^1$ determines the normal arcs of $S\cap \t^2$. For
any tetrahedron $t$ of $\t$, the components of $S\cap t$, being discs,
are determined by $S\cap \d t$, up to isotopy fixed on $\d t$. Thus we
obtain the following lemma. 
\begin{lemma}\label{lem:uniquenormal}
  A pre-normal surface $S\subset M$ is 
  determined by $S\cap \t^1$, up to isotopy mod $\t^2$.\qed
\end{lemma}

\begin{figure}[htbp]
  \begin{center}
    \leavevmode
    \epsfbox{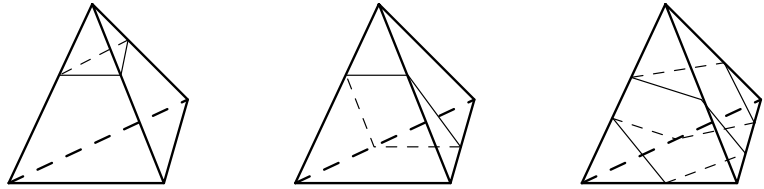}
    \caption{A triangle, a square and an octagon}
    \label{fig:pieces}
  \end{center}
\end{figure}

\begin{definition}\label{def:knormal}
  Let $S\subset M$ be a pre-normal surface and let $k$ be a natural
  number. If for any connected component $C$ of $S\setminus \t^2$ and
  any edge $e$ of $\t$ holds $\#(\d C\cap e)\le k$, then $S$ is
  \begriff{$k$--normal}.
\end{definition}
We are mostly interested in 1-- and 2--normal surfaces. 
Let $S$ be a 2--normal surface and let $t$ be a tetrahedron of $\t$. Then the  
components of $S\cap t$ are copies of triangles, squares and
octagons, as in Figure~\ref{fig:pieces}. 
For any tetrahedron $t$, there are 10 possible types of components
of $S\cap t$: 
four triangles (one for each vertex of $t$), three squares (one for
each pair of opposite edges of $t$), and three octagons. Thus there
are $10\,t(\t)$ possible types of components of $S\setminus
\t^2$. 
Up to isotopy mod $\t^2$, the set $S\setminus \t^2$ is 
described by the vector $\x(S)$ of $10\,t(\t)$ non-negative integers
that indicates the number of copies of the different types of discs
occuring in $S\setminus\t^2$. Note that the 1--normal surfaces are
formed by triangles and squares only. 

We will describe the non-negative integer vectors that represent
2--normal surfaces. Let $S\subset M$ 
be a 2--normal surface and let $x_{t,1},\dots,x_{t,6}$ be the
components of $\x(S)$ that correspond to the squares and octagons in
some tetrahedron $t$. It is impossible that in $S\cap t$ occur two
different types of squares or octagons, since two different squares or
octagons would yield a self-intersection of $S$. Thus  all but
at most one of $x_{t,1},\dots, x_{t,6}$ vanish for any 
$t$. This property of $\x(S)$ is called \begriff{compatibility
  condition}.  

Let $\gamma$ be a normal arc in a 2--simplex $\sigma$ of $\t$ and
$t_1,t_2$ be the two tetrahedra that meet at $\sigma$. In both $t_1$
and $t_2$ there are one triangle, one square and two octagons that
contain a copy of $\gamma$ in its boundary. Moreover, each of them
contains \emph{exactly one} copy of $\gamma$. Let $x_{t_i,1},\dots,
x_{t_i,4}$ be the components of $\x(S)$ that correspond to these types
of discs in $t_i$, where $i=1,2$. 
Since $\d S=\emptyset$, the number of components of
$S\cap t_1$ containing a copy of $\gamma$ equals the number of
components of $S\cap t_2$ containing  a copy of $\gamma$. That is to 
say $x_{t_1,1}+\dots +x_{t_1,4} = x_{t_2,1}+\dots+ x_{t_2,4}$. Thus
$\x(S)$ satisfies a system of linear Diophantine
equations, with one equation for each type of normal arcs. This 
property of $\x(S)$ is called \begriff{matching condition}. The next
claim states that the compatibility and the matching conditions
characterize the vectors that represent 2--normal surfaces. A proof can
be found in~\cite{hemion}, Chapter~9. 
\begin{proposition}\label{prop:vectors}
 Let $\x$ be a vector of $10\,t(\t)$ non-negative integers that
 satisfies  both the compatibility and the matching conditions. Then
 there is a  2--normal surface $S\subset M$ with
 $\x(S)=\x$.\qed 
\end{proposition}

Two 2--normal surfaces $S_1,S_2$ are called \begriff{compatible} if the 
vector $\x(S_1)+\x(S_2)$ satisfies the compatibility condition. It
always satisfies the matching condition. Thus if $S_1$ and $S_2$ are
compatible, then there is a 2--normal surface $S$ with
$\x(S)=\x(S_1)+\x(S_2)$, and we denote $S=S_1+S_2$. 
Conversely, let $S$ be a 2--normal surface, and assume that there are
non-negative integer vectors $\x_1,\x_2$ that both satisfy the
matching condition, with $\x(S)=\x_1+\x_2$. Then both $\x_1$ and
$\x_2$ satisfy the compatibility condition. Thus there are 2--normal
surfaces $S_1,S_2$ with $S=S_1+S_2$. 
The Euler characteristic is additive, i.e., $\chi(S_1+S_2)=\chi(S_1) +
\chi(S_2)$, see~\cite{hemion}. 
\begin{remark}
  The addition of 2--normal surfaces extends to an addition on the set
  of pre-normal surfaces as follows. If $S_1, S_2\subset M$ are
  pre-normal 
  surfaces, then $S_1+S_2$ is the pre-normal surface that is
  determined by $\t^1\cap (S_1\cup S_2)$. The addition yields a
  semi-group structure on the set of pre-normal surfaces. This
  semi-group is isomorphic to the semi-group of integer points in a
  certain rational convex cone that is associated to $\t$. 
  The Euler characteristic is \emph{not} additive with respect to the 
  addition of pre-normal surfaces. 
\end{remark}

\section{Fundamental surfaces}
\label{sec:fundamental}

We use the notations of the previous section.
The power of the theory of 2--normal surfaces is based on the following
two finiteness results.  
\begin{proposition}\label{lem:kneserhaken}
  Let $S\subset M$ be a 2--normal surface comprising more than $10 \,t(\t)$
  two-sided connected components. Then two connected components of $S$ 
  are isotopic mod $\t^2$.\qed 
\end{proposition}
This is proven in~\cite{haken2}, Lemma~4, for 1--normal surfaces. The 
proof easily extends to 2--normal surfaces. 
\begin{theorem}\label{thm:fundamental}
  Let $N\subset M\setminus U(\t^0)$ be a sub--3--manifold whose boundary
  is a 1--normal   surface. 
  There is a system $F_1,\dots, F_q\subset N$ of 2--normal 
  surfaces such that 
  $$\|F_i\| < \|\d N\|\cdot 2^{18 \,t(\t)}$$ for
  $i=1,\dots,q$,
  and any 2--normal surface $F\subset N$ can be written
  as a sum $F=\sum_{i=1}^q k_iF_i$ with non-negative integers
  $k_1,\dots,k_q$.  
\end{theorem}
The surfaces $F_1,\dots, F_q$ are called
\begriff{fundamental}. 
Theorem~\ref{thm:fundamental} is a generalization of a result
of~\cite{hasslagarias} that concerns the case $N=M\setminus
U(\t^0)$.  

The rest of this section is devoted to the proof of
Theorem~\ref{thm:fundamental}. 
The idea is to define a system of linear Diophantine equations
(\emph{matching equations}) whose non-negative solutions correspond to
2--normal surfaces in $N$. 
The fundamental surfaces $F_1,\dots, F_q$ correspond to the Hilbert base
vectors of the equation system, and the bound for $\|F_i\|$ is a
consequence of estimates for the norm of Hilbert base vectors.
Note that in an earlier version of this
paper~\cite{king2}, we proved
Theorem~\ref{thm:fundamental} in essentially the same way, but using handle
decompositions of 3--manifolds rather than triangulations. 

\begin{definition}
  A \begriff{region} of $N$ is a component $R$ of $N\cap
  t$, for a 
  closed tetrahedron $t$ of $\t$. 
  If $\d R\cap \d N$ consists of two copies of one normal triangle or
  normal square then $R$ is a \begriff{parallelity region}.
\end{definition}
\begin{definition}
  The \begriff{class} of a normal triangle, square or
  octagon in $N$ is 
  its equivalence class with respect to isotopies mod $\t^2$ with
  support in $U(N)$.
\end{definition}

Let $t$ be a closed tetrahedron of $\t$, and let $R\subset t$ be a
region of $N$. One verifies that if $R$ is not a parallelity
region then $\d R\cap \d N$ either consists of four normal triangles
(``type~I'')
or of two normal triangles and one normal square (``type~II''). 
If $R$ is of type~I, then $R$ is isotopic mod $\t^2$ to $t\setminus
U(\t^0)$, and any other region of $N$ in $t$ is a parallelity
region. As in the previous section, $R$ contains four classes of normal
triangles, three classes of normal squares and three classes of normal
octagons.  
If $R$ is of type~II, then $t$ contains at most one
other region of $N$ that is not a parallelity region, that is then
also of type~II.
A normal square or octagon in $t$ that is not isotopic mod $\t^2$ to a
component of $\d R\cap \d N$ intersects $\d R$. Thus  $R$ contains two
classes of normal triangles and one class of normal squares.

Let $m(N)$ be the number of classes of normal triangles, squares and
octagons in regions of $N$ of types~I and~II.
If $N$ has $k$ regions of type~I, then $N$ has $\le 2(t(\t)-k)$ regions of 
type~II, thus $m(N) \le 10k + 6(t(\t)-k)\le 10 \,t(\t)$.
Let $\overline m(N)$ be the number of parallelity regions of $N$. It
is easy to see that $\overline m(N) \le \frac 12\, \#(\d N\setminus \t^2) 
\le \frac 16\, \| \d N\| \cdot t(\t)$.

Any 2--normal surface $F\subset N$ is determined up to isotopy mod
$\t^2$ with support in $U(N)$ by the vector $\overline \x_N (F)$ of
$m(N)+ \overline m(N)$ non-negative integers that count the number of
components of $F\setminus \t^2$ in each class of normal triangles,
squares and octagons.  
Let $\gamma_1,\gamma_2\subset \t^2$ be normal arcs, and let
$R_1,R_2$ be two regions of $N$ with $\gamma_1\subset \d
R_1$ and $\gamma_2\subset \d R_2$. 
For $i=1,2$, let $x_{i,1},\dots, x_{i,m_i}$ be the components  of
$\overline \x_N(F)$ that correspond to classes of normal triangles,
squares and octagons in $R_i$ that contain $\gamma_i$ in its boundary. 
  If $x_{1,1} + \dots + x_{1,m_1} = x_{2,1} + \dots + x_{2,m_2}$
  then we say that $\overline \x_N (F)$ satisfies the \begriff{matching equation}
  associated to $(\gamma_1,R_1;\gamma_2,R_2)$.  

For $i=1,2$, $R_i$ contains one class of normal triangles that contain
a copy of $\gamma_i$ in its boundary. If $R_i$ is not a parallelity
region, then $R_i$ contains one class of normal squares that contain a
copy of $\gamma_i$ in its boundary. If $K_i$ is of type~I, then $K_i$
additionally contains two classes of normal octagons containing a copy
of $\gamma_i$ in its boundary. Thus if $R_i$ is a parallelity region
then $m_i=1$, if it is of type~I then $m_i=4$, and if it is of type~II 
then $m_i=2$.

For any 2--normal surface $F\subset N$, let
$\x_N(F)\in \mathbb Z_{\ge 0}^{m(N)}$ be the vector that collects the
components of $\overline\x_N(F)$ corresponding to the classes of normal
triangles, squares and octagons in regions of $N$ of types~I and~II. 
As in the previous section, the vector $\x_N(F)$ (resp.\
$\overline\x_N(F)$) satisfies a \emph{compatibility condition}, i.e., 
for any region $R$ of $N$ vanish all but at most one components of
$\x_N(F)$ (resp. $\overline \x_N(F)$) corresponding to 
classes of squares and octagons in $R$. 

\begin{lemma}\label{lem:vectors2}
  Suppose that any component of $N$ contains a region that is not a
  parallelity region.  
  There is a system of matching equations concerning only regions of
  $N$ of types~I and~II, such that a vector $\x\in \mathbb Z_{\ge
    0}^{m(N)}$ satisfies these equations and the compatibility
  condition if and only if there is a 2--normal surface $F\subset N$
  with $\x_N(F) = \x$. The surface $F$ is determined by $\x_N(F)$, up
  to isotopy in $N$ mod $\t^2$.
\end{lemma}
\begin{proof}
  Let $\gamma\subset N\cap \t^2$ be a normal arc. Let $R_1,R_2$ be the
  two regions of $N$ that contain $\gamma$.
  Let $F\subset N$ be a 2--normal surface. Since $\d F=\emptyset$, the
  number of components of $F\cap R_1$ containing $\gamma$ and the number
  of components of $F\cap R_2$ containing $\gamma$ coincide.
  Thus $\overline\x_N(F)$ satisfies the matching equation associated to
  $(\gamma,R_1;\gamma,R_2)$. We refer to these equations as
  $N$--matching equations. We will transform the system of $N$--matching
  equations 
  by eliminating the components of $\overline\x_N(F)$ that do not
  belong to $\x_N(F)$.

  Let $\gamma_1,\gamma_2\subset \t^2$ be normal arcs, and let
  $R_1,R_2$ be two different regions of $N$ with $\gamma_1\subset \d
  R_1$ and $\gamma_2\subset \d R_2$. 
  Assume that $R_1$ is a parallelity region of $N$. Then $m_1=1$, thus
  the matching equation associated to $(\gamma_1,R_1;\gamma_2,R_2)$ is
  of the form $ x_{1,1} = x_{2,1} + \dots + x_{2,m_2} $. Hence we can
  eliminate $x_{1,1}$ in the $N$--matching equations. 
  For any region $R_3$ of $N$ and any normal arc $\gamma_3\subset \d
  R_3$, the elimination transforms the matching equation associated to
  $(\gamma_1,R_1;\gamma_3,R_3)$ into the matching equation associated to 
  $(\gamma_2,R_2;\gamma_3,R_3)$. 
  We iterate the elimination process. Since any component of $N$
  contains a region that is not a parallelity region, we eventually
  transform the system of $N$--matching equations to a system
  $\mathfrak A$ of matching equations that concern only regions of $N$
  of types~I and~II. 

  Let $\x\in \mathbb Z_{\ge 0}^{m(N)}$ be a solution of $\mathfrak
  A\cdot \x = 0$. By the elimination process, there is a unique
  extension of $\x$ to a solution
  $\overline \x$ of the $N$--matching equations. If $\x$ satisfies the
  compatibility condition then so does $\overline \x$, since a
  parallelity region contains at most one class of normal squares. 
  Now the lemma follows by 
  Proposition~\ref{prop:vectors}, that is proven in~\cite{hemion}.
\end{proof}

\paragraph*{Proof of Theorem~\ref{thm:fundamental}}
  It is easy to verify that if $R$ is a parallelity region then there
  is only one class of 2--normal pieces in $R$. 
  If a component $N_1$ of $N$ is a union of parallelity regions, then
  $N_1$ is a regular neighbourhood of a 1--normal surface $F_1\subset
  N_1$, that has a connected non-empty intersection with each region
  of $N_1$. Any pre-normal surface in $N_1$ is a multiple of  
  $F_1$ (thus, is 1-normal), see~\cite{haken1}. 
  We have $\|F_1\| = \frac 12\, \|\d N_1\|$. 
  Thus by now we can suppose that any component of $N$ contains a
  region that is not a parallelity region.

  By Lemma~\ref{lem:vectors2}, the $\x$--vectors of 2--normal surfaces
  in $N$ satisfy a system of linear equations $\mathfrak
  A\cdot\x=0$. By the following well known result on Integer Programming
  (see~\cite{schrijver}), the non-negative integer solutions of such a
  system are additively generated by a finite set of solutions. 
\begin{lemma}\label{lem:hilbert}
  Let $\mathfrak A = (a_{ij})$ be an integer $(n\times
  m)$--matrix. 
  Set $$K=\left( \max_{i=1,\dots,n} \sum_{j=1}^m
                                  a_{ij}^2\right)^{1/2}.$$ 
  There is a set $\{\x_1,\dots \x_p\}$ of non-negative integer
  vectors such that $\mathfrak A\cdot \x_i = 0$ for any $i=1,\dots,p$,
  the components of $\x_i$ are bounded from above by $mK^m$, 
  and any non-negative 
  integer solution $\x$ of $\mathfrak A\cdot \x = 0$ can
  be written as a sum $\x = \sum k_i\x_i$ with non-negative integers
  $k_1,\dots,k_p$.\qed
\end{lemma}
The set $\{\x_1,\dots \x_p\}$ is called
\begriff{Hilbert base}\index{Hilbert base} for $\mathfrak A$, if $p$
is minimal.  

  As in the previous section, if $F\subset N$ is a 2--normal surface
  and $\x_{N}(F)$ is a sum of two non-negative integer solutions of
  $\mathfrak A\cdot\x=0$ then there are 2--normal surfaces 
  $F',F''\subset N$ with $F=F'+F''$. 
  Thus the surfaces $F_1,\dots, F_q\subset N$ that correspond to
  Hilbert base vectors satisfying the compatibility condition 
  additively generate the set of all 2--normal surfaces in $N$.

  It remains to bound $\|F_i\|$, for $i=1,\dots, q$. Since $F_i$ is
  2--normal and any edge of $\t$ is of degree $\ge 3$, 
  we have $\|F_i\|\le \frac 83\, \#(F_i\setminus \t^2)$.
  By the elimination process in the proof of Lemma~\ref{lem:vectors2},
  any component of $\overline \x_N(F_i)$ that corresponds to a
  parallelity region of $N$ is a sum of at most four components of
  $\x_N(F_i)$. 
  By the bound for the components of $\x_N(F_i)$ in
  Lemma~\ref{lem:hilbert} (with $m=m(N)$ and $K^2=8$) and our 
  bounds for $m(N)$ and $\overline m(N)$, we obtain 
  \begin{eqnarray*}
   \|F_i\| &\le& \frac 83 \cdot \left(m(N) + 4\,\overline m(N)\right) \cdot
                               \left(m(N)\cdot 2^{\frac 32 m(N)}\right) \\
           &\le& \frac 83\cdot \left( 10 \,t(\t) + \frac 23\, \|\d N\|\, t(\t)\right)\cdot
           10\,t(\t)\cdot 2^{15 \,t(\t)} \\
           &<& \left(300 + 20\, \| \d N\| \right)\cdot t(\t)^2\cdot 2^{15 \,t(\t)} .
  \end{eqnarray*}
  Using $t(\t)\ge 5$ and $\|\d N\|>0$, we obtain $\|F_i\|< \|\d N\|\cdot 
  2^{18 \,t(\t)}$.\qed


\section{Maximal systems of  1--normal spheres}
\label{sec:constrN}

Let $\t$ be a triangulation of a closed orientable 3--manifold $M$. 
By Proposition~\ref{lem:kneserhaken}, there is a system $\Sigma\subset
M$ of $\le 10\,t(\t)$ pairwise disjoint 1--normal spheres, such that any
1--normal sphere in $M\setminus \Sigma$ is isotopic mod $\t^2$ to a
component of $\Sigma$. We call such a system \begriff{maximal}. It is
not obvious how to {construct} $\Sigma$, in particular how to
estimate $\|\Sigma\|$ in terms of $t(\t)$. 
This section is devoted to a solution of this problem.

\begin{construction}\label{con:Sigma}
  Set $\Sigma_1=\d U(\t^0)$ and $N_1=M\setminus U(\t^0)$. 
  Let $i\ge 1$. If there is a 1--normal fundamental projective plane 
  $P_i\subset N_i$ then set $\Sigma_{i+1}=\Sigma_i\cup 2P_i$ and
  $N_{i+1} = N_i\setminus U(P_i)$.
  Otherwise, if there is a 1--normal fundamental sphere $S_i\subset
  N_i$ that is not isotopic mod $\t^2$ to a component of $\Sigma_i$,  
  then set $\Sigma_{i+1}=\Sigma_i\cup S_i$ and $N_{i+1}=N_i\setminus
  U(S_i)$. 
  Otherwise, set $\Sigma=\Sigma_i$. 
\end{construction}
Since $M$ is orientable, a projective plane $P_i$ is
one-sided and $2P_i$ is a sphere.
By Proposition~\ref{lem:kneserhaken} and since
embedded spheres are two-sided in $M$, the iteration
stops for some $i < 10 \,t(\t)$. 

\begin{lemma}\label{lem:boundsigma}
  $\|\Sigma\| < 2^{185 \,t(\t)^2}$.
\end{lemma}
\begin{proof}
  In Construction~\ref{con:Sigma}, we have  
  \begin{eqnarray*}
    \|\Sigma_{i+1}\| &<& \|\Sigma_i\| + 2 \|\Sigma_i\|\cdot 2^{18 \,t(\t)}
        \\ 
      &<& \|\Sigma_i\|\cdot 2^{18 \,t(\t) +2}
  \end{eqnarray*}
  by Theorem~\ref{thm:fundamental} .
  The iteration stops after $< 10 \,t(\t)$ steps, thus $$\|\Sigma\| 
  < \|\Sigma_1\| \cdot 2^{180 \,t(\t)^2 + 20 \,t(\t)} \le \|\Sigma_1\| \cdot
  2^{184 \,t(\t)^2},$$ using $t(\t)\ge 5$.  
  Since $\| \d U(\t^0)\|$ equals twice the number of edges of $\t$, we
  have $\|\Sigma_1\| \le 4 \,t(\t)$, and the lemma
  follows.
\end{proof}

\begin{lemma}\label{lem:keinS1}
  $\Sigma$ is maximal.
\end{lemma}
\begin{proof}
  It is to show that any 1--normal sphere $S\subset M\setminus
  U(\Sigma)$ is isotopic mod $\t^2$ to a component of $\Sigma$.
  Let $N$ be the component of $M\setminus U(\Sigma)$ that contains
  $S$.
  If $N$ contains a 1--normal fundamental projective plane $P$, then
  $N=U(P)$ by Construction~\ref{con:Sigma}. Thus $S=2P=\d N$, which is
  isotopic mod $\t^2$ to a component of $\Sigma$. Hence we can assume
  that $N$ does not contain a 1--normal fundamental projective plane.

  We express $S$ as a sum $\sum_{i=1}^q k_iF_i$ of fundamental
  surfaces in $N$. 
  Since $\chi(S)=2$ and the Euler characteristic is additive, one of
  the fundamental surfaces in the sum, say, $F_1$ with $k_1>0$, has
  positive Euler characteristic. It is not a projective plane by the
  preceding paragraph, thus it is a sphere. 
  By construction of $\Sigma$, the sphere $F_1$ is isotopic mod $\t^2$
  to a component of $\Sigma$, thus it is parallel to a component of
  $\d N$. 
  Hence $F_1$ is disjoint to any 1--normal surface in $N$, up to
  isotopy mod $\t^2$. Thus $S$ is the disjoint union of $k_1F_1$ and
  $\sum_{i=2}^q k_i F_i$. Since $S$ is connected, it follows
  $S=F_1$. Thus $S$ is isotopic mod $\t^2$ to a component of
  $\Sigma$. 
\end{proof}

We will extend $\Sigma$ to a system $\tilde\Sigma$ of 2--normal
spheres. To define $\tilde\Sigma$, we need a lemma about 2--normal
spheres in the complement of $\Sigma$.
\begin{lemma}\label{lem:boundS2}
  Let $N$ be a component of $M\setminus U(\Sigma)$. Assume that
  there is a 2--normal sphere in $N$ with exactly one
  octagon. Then there is a 2--normal fundamental sphere $F\subset N$
  with exactly one octagon and $\|F\| < 2^{189 \,t(\t)^2}$.
\end{lemma}
\begin{proof}
  Let $S\subset N$ be a 2--normal sphere with exactly one octagon.
  If $N$ contains a 1--normal fundamental projective plane $P$, then
  $N=U(P)$ by Construction~\ref{con:Sigma}, and any pre-normal surface
  in $N$ is a multiple of $P$, i.e., is 1--normal.
  Thus since $S\subset N$ is not 1--normal, there is no 
  1--normal fundamental projective plane in $N$.

  We write $S$ as a sum of 2--normal fundamental surfaces in $N$. 
  Since $S$ has exactly one octagon, exactly one summand is not
  1--normal. Since any projective plane in the sum is not 1--normal by
  the preceding paragraph, at most one summand is a projective plane.
  Since $\chi(S)=2$ and the Euler characteristic is additive, it
  follows that one of the fundamental surfaces in the sum is a sphere
  $F$. 

  Assume that $F$ is 1--normal, i.e.,  $S\not=F$. 
  The construction of $\Sigma$ implies that $F$ is isotopic mod $\t^2$
  to a component of $\d N$. Thus it is disjoint to any
  2--normal surface in $N$. 
  Therefore $S$ is a disjoint union of a multiple of $F$ and of a
  2--normal surface with exactly one octagon, which is a
  contradiction since $S$ is connected. Hence $F$ contains the
  octagon of $S$. 
  We have $\|F\| <  \|\Sigma\|\cdot 2^{18 \,t(\t)}$ by
  Theorem~\ref{thm:fundamental}. 
  The claim follows with Lemma~\ref{lem:boundsigma} and $t(\t)\ge
  5$.  
\end{proof}

The preceding lemma assures that the following construction works.
\begin{construction}\label{con:tildeSigma}
  For any connected component $N$ of $M\setminus U(\Sigma)$ that
  contains a 2--normal sphere with exactly one octagon, choose a
  2--normal sphere $F_N\subset N$ with exactly one octagon and  $\|F\|
  < 2^{189 \,t(\t)^2}$. Set
     $$\tilde\Sigma = \Sigma\cup\bigcup_N F_N.$$ 
\end{construction}
Since $\#(\tilde\Sigma) \le 10\,t(\t)$ by
Proposition~\ref{lem:kneserhaken}, it follows $\|\tilde\Sigma\| <
10\,t(\t)\cdot 2^{189 \,t(\t)^2} < 2^{190 \,t(\t)^2}$.


\section{Almost $k$--normal surfaces and split equivalence}
\label{sec:almost}

We shall need a generalization of the
notion of $k$--normal surfaces. Let $M$ be a closed connected
orientable 3--manifold with a triangulation $\t$. 

\begin{definition}\label{def:akn}
  A closed embedded surface $S\subset M$ transversal to $\t^2$ is
  \begriff{almost $k$--normal}, if 
  \begin{enumerate}
  \item $S\cap \t^2$ is a union of normal arcs and of circles in
    $\t^2\setminus\t^1$, and
  \item for any tetrahedron $t$ of $\t$, any edge $e$ of $t$ and any
    component $\zeta$ of $S\cap \d t$ holds $\#(\zeta\cap e)\le k$. 
  \end{enumerate}
\end{definition}

Our definition is similar to Matveev's one in~\cite{matv}. Note
that there is a related but different definition of \lq\lq almost
normal\rq\rq\ surfaces due to Rubinstein~\cite{rubinstein}. 
Any surface in $M$ disjoint to $\t^1$ is almost 1--normal. Any almost
$k$--normal surface that meets $\t^1$ can be seen as a
$k$--normal surface with several disjoint small tubes 
attached in $M\setminus \t^1$, see~\cite{matv}. The tubes can be
nested. Of course there are many ways to add tubes to a  $k$--normal
surface. We shall develop tools to deal with this ambiguity. 

Let $S\subset M$ be an almost $k$--normal surface. By definition, the
connected components of $S\cap \t^2$ that meet $\t^1$ are formed by
normal arcs. Thus these components define a pre-normal surface
$S^\times$, that is obviously $k$--normal. It is determined by $S\cap
\t^1$, according 
to Lemma~\ref{lem:uniquenormal}. 
A disc $D\subset M\setminus \t^1$ with $\d D\subset S$  is called a
\begriff{splitting disc} for $S$. One obtains
$S^\times$ by splitting $S$ along splitting discs for $S$ that are
disjoint to $\t^2$, and isotopy mod $\t^1$.

If two almost $k$--normal surfaces $S_1, S_2$ satisfy
$S_1^\times=S_2^\times$, then $S_1$ and $S_2$ differ only by the
choice of tubes. This gives rise to the following equivalence
relation. 
\begin{definition}
  Two embedded  surfaces $S_1,S_2\subset M$ transversal to
  $\t^2$ are \begriff{split equivalent} if 
  $S_1\cap\t^1=S_2\cap\t^1$ (up to isotopy  mod $\t^2$). 
\end{definition}
If two almost $k$--normal surfaces $S_1,S_2\subset M$ are split
equivalent, then $S_1^\times =S_2^\times$, up
to isotopy mod $\t^2$. In particular, two $k$--normal surfaces are
split equivalent if and only if they are isotopic mod $\t^2$.

\begin{definition}
  If $S\subset M$ is an almost $k$--normal surface and $S^\times$ is
  the disjoint union of $k$--normal surfaces
  $S_1,\dots, S_n$, then we call $S$ a \begriff{tube sum} of
  $S_1,\dots,S_n$. We denote the set of all tube sums of
  $S_1,\dots,S_n$ by $S_1\circ\dots\circ S_n$.
\end{definition}
\begin{definition}
  Let  $S=S_1\cup\dots\cup S_n\subset M$ be a surface transversal to
  $\t^2$ with $n$
  connected components, and let $\Gamma\subset M\setminus \t^1$ be a
  system of disjoint simple arcs with $\Gamma\cap S= \d \Gamma$. For
  any arc $\gamma$ in $\Gamma$, one component  of $\d U(\gamma)
  \setminus S$ is an annulus $A_\gamma$. 
  The surface $$S^\Gamma = (S\setminus U(\Gamma)) \cup
  \bigcup_{\gamma\subset 
    \Gamma} A_\gamma$$ is called the 
  \begriff{tube sum of $S_1,\dots, S_n$ along $\Gamma$}. 
\end{definition}
If $S_1,\dots,S_n$ are $k$--normal, then $S^\Gamma\in S_1\circ
\dots\circ S_n$.  

We recall the concept of impermeable surfaces, that is central in the
study of almost 2--normal surfaces (see~\cite{thompson},\cite{matv}). 
Fix a vertex $x_0\in\t^0$.
Let $S\subset M$ be a connected embedded surface transversal to
$\t$. If $S$ splits $M$ into two pieces, then let $B^+(S)$ denote the 
closure of the  component of $M\setminus S$ that contains $x_0$,
and let $B^-(S)$ denote the closure of the other component. We do not
include $x_0$ in the notation ``$B^+(S)$'', since in our  
applications the choice of $x_0$ plays no essential role.

\begin{definition}
  Let $S\subset M$ be a connected embedded surface transversal to
  $\t^2$. Let 
  $\alpha\subset\t^1\setminus\t^0$ and $\beta\subset S$ be embedded
  arcs with $\d\alpha=\d \beta$. A closed embedded disc $D\subset M$
  is a \begriff{compressing disc} for $S$ with {string} $\alpha$ and
  base $\beta$, if $\d D=\alpha\cup\beta$ and $D\cap \t^1=\alpha$. If, 
  moreover, $D\cap S=\beta$, then we call $D$ a \begriff{bond} of
  $S$. 

  Let $S\subset M$ be a connected embedded surface that splits
  $M$ and let $D$ be a compressing disc for $S$ with string
  $\alpha$. If the germ of $\alpha$ in $\d\alpha$ is contained in
  $B^+(S)$ (resp.\ $B^-(S)$), then  $D$ is \begriff{upper} (resp.\
  \begriff{lower}).  
  Let $D_1, D_2$ be upper and lower compressing discs for $S$
  with strings $\alpha_1,\alpha_2$. If $D_1\subset D_2$ or $D_2\subset 
  D_1$, then $D_1$ and $D_2$ are \begriff{nested}. If $D_1\cap
  D_2\subset \d\alpha_1\cap \d\alpha_2$, then $D_1$ and $D_2$ are
  \begriff{independent} from each other.
\end{definition}
Upper and lower compressing discs that are independent from each
other meet in at most one point. 

\begin{definition}
  Let $S\subset M$ be a connected embedded surface that is transversal
  to $\t^2$ and splits $M$. 
  If $S$ has both upper and lower bonds, but no pair of
  nested or independent upper and lower compressing discs, then $S$ is
  \begriff{impermeable}. 
\end{definition}
Note that the impermeability of $S$ does not change under an isotopy
of $S$ mod $\t^1$. The next two claims state a close relationship
between impermeable surfaces and (almost) 2--normal surfaces. 
By an octagon of an almost 2--normal surface $S\subset M$ in a tetrahedron 
$t$, we mean a circle in $S\cap \d t$ formed by eight normal arcs. This
corresponds to an octagon of $S^\times$ in the sense of
Figure~\ref{fig:pieces}. 
\begin{proposition}\label{prop:impermeabel}
  Any impermeable surface in $M$ is isotopic mod $\t^1$ to an almost
  2--normal surface with exactly one octagon.  
\end{proposition}
\begin{proposition}\label{prop:2n-imperm}
  A connected 2--normal surface that splits $M$ and contains exactly one
  octagon is impermeable.
\end{proposition}
We shall need these statements later.
As the author found only parts of the proofs in 
the literature (see~\cite{thompson},\cite{matv}), he includes proofs
in Section~\ref{sec:imperm-a2}. 

We end this section with the
definition of $\t^1$--Morse embeddings and with 
the notion of thin position. Let $S$ be a closed 2--manifold and let
$H\co S\times I\to M$ be a tame embedding. For $\xi\in I$, set
$H_\xi=H(S\times \xi)$.  
\begin{definition}
  An element $\xi\in I$ is a \begriff{critical parameter} of $H$ and a
  point $x\in H_\xi$ is a \begriff{critical point} of $H$ with respect
  to $\t^1$, if $x$ is a vertex of $\t$ or $x$ is a point of tangency
  of $H_\xi$ to $\t^1$. 
\end{definition}
\begin{definition}\label{def:regular}
  We call $H$ a \begriff{$\t^1$--Morse embedding}, if it has finitely many
  critical parameters, to any critical parameter belongs exactly one
  critical point, and for any critical point $x\in \t^1\setminus
  \t^0$ corresponding to a critical parameter $\xi$, one component of
  $U(x)\setminus H_\xi$ is disjoint to $\t^1$. The number of critical
  points with respect to $\t^1$ of a $\t^1$--Morse embedding $H$ is denoted
  by $c(H,\t^1)$. 
\end{definition}
The last condition in the definition of $\t^1$--Morse embeddings means that
any critical point of $H$ is a vertex of $\t$ or a local maximum
resp.\ minimum of an edge of $\t$.
\begin{definition}
  Let $F$ be a closed surface, let $J\co F\times I\nach M$ be a $\t^1$--Morse
  embedding, and let $\xi_1,\dots,\xi_r\in I$ be 
  the critical parameters of $J$ with respect to $\t^1$. The
  \begriff{complexity} $\kappa(J)$ of $J$ is defined as 
  $$\kappa(J) = \#\left(\t^1\setminus \left( \bigcup_{i=1}^r
      J_{\xi_i}\right)\right).$$  
\end{definition}
If $\kappa(J)$ is minimal among all $\t^1$--Morse
embeddings with the property $\t^1\subset J(F\times I)$, then $J$ is
said to be in \begriff{thin position} with respect to $\t^1$.
This notion was introduced for foliations of 3-manifolds by 
Gabai~\cite{gabai}, was applied by Thompson~\cite{thompson} for her 
recognition algorithm of $S^3$, and was also used in the study of
Heegaard surfaces by Scharlemann and Thompson~\cite{scharlthompson}.  

If $J(F\times\xi)$ splits $M$ and has a pair of nested or
independent upper and lower compressing discs $D_1,D_2$, then an
isotopy of $J$ along $D_1\cup D_2$ decreases $\kappa(J)$,
see~\cite{matv}, \cite{thompson}. We obtain the following
claim. 
\begin{lemma}
  Let $J\co F\times I\to M$ be a $\t^1$--Morse embedding in thin position
  and let $\xi\in I$ be a non-critical parameter of $J$. If $J(F\times
  \xi)$ has both upper and lower bonds, then $J(F\times \xi)$ is
  impermeable. \qed  
\end{lemma}


\section{Compressing and splitting discs}
\label{sec:1-2-normal}

Let $M$ be a closed connected 3--manifold with a triangulation $\t$. In
the lemmas that we prove in this section, we state technical conditions
for the existence of compressing and splitting discs for a surface.

\begin{lemma}\label{lem:keinelcd}
  Let $S_1,\dots,S_n\subset M$ be embedded surfaces transversal to
  $\t^2$ and let $S$ be
  the tube sum of $S_1,\dots, S_n$ along a system $\Gamma\subset
  M\setminus \t^1$ of arcs. 
  Assume that $S$ splits $M$, and $\Gamma\subset B^-(S)$. If none of
  $S_1,\dots,S_n$ has a lower compressing disc, 
  then $S$ has no lower compressing disc.
\end{lemma}
\begin{proof}
  Set $\Sigma = S_1\cup\dots\cup S_n$.
  Let $D\subset M$ be a lower compressing disc for $S$. 
  One can assume that a collar of $\d D\cap S$ in $D$ is contained in
  $B^-(S)$. 
  Then, since by hypothesis $U(\Gamma)\cap \Sigma\subset B^-(S)$,
  any point in $\d D\cap U(\Gamma)\cap \Sigma$ is endpoint of an arc
  in $D\cap \Sigma$. Therefore there is a sub-disc
  $D'\subset D$, bounded by parts of $\d D$ and of arcs in
  $D\cap\Sigma$, that is a lower compressing disc for one of
  $S_1,\dots,S_n$. 
\end{proof}

\begin{lemma}\label{lem:makeindep}
  Let $S\subset M$ be a surface transversal to $\t^2$ with upper and
  lower compressing discs 
  $D_1$, $D_2$ such that $\d (D_1\cap D_2)\subset \d D_2\cap
  S$. Assume either that $(\d D_1)\cap D_2\subset \t^1$ or 
  that there is a splitting disc $D_m$ for $S$ such that $D_1\cap D_m=
  \d D_1\cap \d D_m=\{x\}$ is a single point and $D_2\cap
  D_m=\emptyset$. Then $S$ has a pair of independent or nested upper and
  lower compressing discs.
\end{lemma}
\begin{proof}
  If $D_1\cap D_2\cap\t^1$ comprises more than a single point then the
  string of $D_2$ is contained in the string of $D_1$. Thus $D_1\cap S$
  contains an arc different from the base of $D_1$, bounding in $D_1$
  a lower compressing disc, that forms with $D_1$ a pair of nested upper
  and lower compressing discs for $S$.

  Assume that a component $\gamma$ of $D_1\cap D_2$ is a circle. Then
  there are discs $D_1'\subset 
  D_1$ and $D_2'\subset D_2$ with $\d D_1'=\d D_2'=\gamma$. 
  Since $\d (D_1\cap D_2)\subset \d D_2$ , $D_2'$ does not contain arcs
  of $D_1\cap D_2$. 
  Thus 
  if we choose $\gamma$ innermost in $D_2$, then $D_1\cap
  D_2'=\gamma$. By cut-and-paste of $D_1$ along $D_2'$, one reduces
  the number of circle components in $D_1\cap D_2$. Therefore we
  assume by now that $D_1\cap D_2$ consists of isolated points in $\d
  D_1\cap \d D_2$ and of arcs that do not meet $\d D_1$.

  Assume that there is a point $y\in (\d D_1\cap \d D_2)\setminus
  \t^1$. Then there is an arc $\gamma \subset \d D_1$ with
  $\d\gamma=\{x,y\}$. Without 
  assumption, let $\gamma\cap D_2=\{y\}$. Let $A$ be the closure of the
  component of $U(\gamma)\setminus (D_1\cup D_2\cup D_m)$ whose
  boundary contains arcs in both $D_2$ and $D_m$. Define $D^*_2=
  ((D_2\cup D_m)\setminus U(\gamma)) \cup A$, that is to say, $D^*_2$ 
  is the connected sum of $D_2$ and $D_m$ along $\gamma$. By
  construction, $(D_1\cap D^*_2)\setminus \d D_1=(D_1\cap
  D_2)\setminus \d D_1$, and 
  $\#(D_1\cap D^*_2) < \#(D_1\cap D_2)$. In that way, we remove all
  points of intersection of $(\d D_1\cap D_2)\setminus
  \t^1$. 
  Thus by now we can assume that $D_1\cap D_2$ consists of
  arcs in $D_1$ that do not meet $\d D_1$, and possibly of a single
  point in $\t^1$. 

  Let $\gamma\subset D_1\cap D_2$ be an outermost arc
  in $D_2$, that is to say, $\gamma\cup \d D_2$ bounds a disc
  $D'\subset D_2\setminus \t^1$ with $D_1\cap D'=\gamma$. We move
  $D_1$ away from $D'$ by an isotopy mod $\t^1$ and obtain a
  compressing disc $D^*_1$ for $S$ with $D^*_1\cap 
  D_2= (D_1\cap D_2)\setminus\gamma$. In that way, we remove all arcs
  of $D_1\cap D_2$ and finally get a pair of independent upper and
  lower compressing discs for $S$. 
\end{proof}

\begin{lemma}\label{lem:1ohnedisc}
  Let $S\subset M$ be an almost 1--normal surface. If $S$ has a
  compressing  disc, then $S$ is isotopic mod $\t^1$ to an
  almost 1--normal surface with a compressing disc contained in a
  single tetrahedron. In particular, $S$ is not 1--normal. 
\end{lemma}
\begin{proof} 
  Let $D$ be a compressing disc for $S$. Choose $S$ and $D$ up to
  isotopy of $S\cup D$ mod $\t^1$ 
  so that $S$ is almost 1--normal and $\#(D\cap\t^2)$ is
  minimal. Choose an innermost component $\gamma\subset(D\cap\t^2)$,
  which is possible as $D\cap \t^2\not=\emptyset$.
  There is a closed tetrahedron $t$ of $\t$ and a
  component $C$ of $D\cap t$ that 
  is a disc, such that $\gamma=C\cap \d t$. Let $\sigma$ be
  the closed 2--simplex of $\t$ that contains $\gamma$. We obtain three
  cases. 
  \begin{enumerate}
  \item Let $\gamma$ be a circle, thus $\d C=\gamma$. Then there
    is a disc $D'\subset\sigma$ with $\d D'=\gamma$ and a ball
    $B\subset t$ with $\d B=C\cup D'$. By an
    isotopy mod $\t^1$ with support in $U(B)$, we move $S\cup
    D$ away from $B$, obtaining a surface $S^*$ with a compressing disc
    $D^*$. If $S^*$ is almost 1--normal, then we obtain a
    contradiction to our choice as 
    $\#(D^*\cap \t^2)<\#(D\cap \t^2)$. 
  \item Let $\gamma$ be an arc with endpoints in a single component
    $c$ of $S\cap\sigma$. Since $S$ has no returns, $\gamma$ is not 
    the string of $D$. We apply to $S\cup D$ an isotopy mod
    $\t^1$ with support in $U(C)$ that moves $C$ into $U(C)\setminus
    t$, and obtain a surface $S^*$ with a compressing disc $D^*$. If
    $S^*$ is almost 1--normal, then we obtain a contradiction to our
    choice as $\#(D^*\cap \t^2)<\#(D\cap \t^2)$. 
  \item Let $\gamma$ be an arc with endpoints in two different
    components $c_1,c_2$ of $S\cap\sigma$. If both $c_1$ and $c_2$ are
    normal arcs, then set $C'=C$, $c'_1=c_1$ and $c'_2=c_2$. If, say,
    $c_1$ is a circle, then we move $S\cup D$ away from $C$ by an isotopy
    mod $\t^1$ with support in $U(C)$. If the resulting surface
    $S^*$ is still almost 1--normal, 
    then we obtain a contradiction to the choice of $D$.  
  \end{enumerate}
  In either case, $S^*$ is not almost 1-normal, i.e., the isotopy introduces a
  return. Therefore there is a component of $C\setminus S$ 
  with closure $C'$ such that $\d C'\cap S$ connects two normal
  arcs $c'_1,c'_2$ of $S\cap\sigma$.  

  Let $\gamma'=C'\cap \sigma$.
  Up to isotopy of $C'$ mod $\t^2$ that
  is fixed on $\d C'\cap S$, we assume that $\gamma'\cap (c'_1\cup
  c'_2)\subset \d \gamma'$. There is an arc $\alpha$ contained in an
  edge of $\sigma$ with $\d\alpha\subset c'_1\cup c'_2$. For
  $i\in\{1,2\}$, there is an arc $\beta_i\subset c'_i$ that connects
  $\alpha\cap c'_i$ with $\gamma'\cap c'_i$. The circle
  $\alpha\cup\beta_1\cup\beta_2\cup\gamma'$ bounds a closed disc
  $D'\subset \sigma$. Eventually $D'\cup C'$ is a compressing disc for
  $S$ contained in a single tetrahedron.   
\end{proof}

\begin{lemma}\label{lem:discparallel}
  Let $S\subset M$ be a 1--normal surface and let $D$ be a splitting
  disc for $S$. Then, $(D,\d D)$ is isotopic in $(M\setminus \t^1,
  S\setminus\t^1)$ to a disc embedded in $S$.
\end{lemma}
\begin{proof}
  We choose $D$ up to isotopy of $(D,\d D)$ in $(M\setminus \t^1,
  S\setminus\t^1)$ so that $(\#((\d D)\cap\t^2), \#(D\cap \t^2))$ is
  minimal in lexicographic order. Assume that $\d D\cap
  \t^2\not=\emptyset$. Then, there is a tetrahedron $t$, a
  2--simplex $\sigma\subset \d t$, a component $K$ of $S\cap t$, and a
  component $\gamma$ of $\d D\cap K$ with $\d \gamma\subset
  \sigma$. Since $S$ is 1--normal, the closure $D'$ of one
  component of $K\setminus \gamma$ is a disc with $\d D'\subset
  \gamma\cup\sigma$. By choosing $\gamma$ innermost in $D$, we can
  assume that $D'\cap \d D = \gamma$. An isotopy of $(D,\d D)$ in
  $(M\setminus \t^1, S\setminus\t^1)$ with support in $U(D')$, moving
  $\d D$ away from $D'$, reduces $\#(\d D\cap \t^2)$, in contradiction
  to our choice. Thus $\d D\cap \t^2=\emptyset$. 

  Now, assume that $D\cap\t^2\not=\emptyset$. Then, there is a
  tetrahedron $t$, a 2--simplex $\sigma\subset \d t$, and a disc
  component $C$ of $D\cap t$, such that $C\cap \sigma=\d C$ is a
  single circle. There is a ball $B\subset t$ bounded by $C$ and a
  disc in $\sigma$. An isotopy of $D$ with support in $U(B)$, moving
  $C$ away from $t$, reduces $\#(D\cap \t^2)$, in contradiction to our
  choice. Thus $D$ is contained in a single tetrahedron $t$. Since $S$
  is 1--normal, $\d D$ bounds a disc $D'$ in $S\cap t$. An isotopy with
  support in $t$ that is constant on $\d D$ moves $D$ to $D'$, which
  yields the lemma.
\end{proof}

\begin{corollary}\label{cor:gibtJ}
  Let $S_0\subset M$ be a 1--normal sphere that splits $M$, and let
  $S\subset B^-(S_0)$  be an almost 1--normal sphere disjoint to $S_0$
  that is split equivalent to $S_0$. Then 
  there is a $\t^1$--Morse embedding $J\co S^2\times I\to M$ with
  $J(S^2\times I) = B^+(S)\cap B^-(S_0)$ and $c(J,\t^1) = 0$.
\end{corollary}
\begin{proof}
  Let $X$ be a graph isomorphic to $S_0\cap \t^2$. Since $S^\times$ is a
  copy of $S_0$, there is an embedding $\varphi\co X\times I \to B^+(S)\cap
  B^-(S_0)$ with  
  $\varphi(X^0\times I) = \varphi(X\times I)\cap \t^1$, $\varphi(X\times
  0) = S_0\cap \t^2 = S_0\cap  \varphi(X\times I)$, and $\varphi(X\times
  1)$ is the union of the normal arcs in $S$. 

  Let $\gamma\subset S\cap \varphi(X\times I)$ be a circle that does not
  meet $\t^1$. 
  Then, $\gamma$ bounds a disc $D\subset \varphi(X\times I)\setminus
  \t^1$. The two components of $S\setminus \gamma$ are discs. One of
  them is disjoint to $\t^1$, since otherwise the disc $D$ would give
  rise to a splitting disc for $S^\times = S_0$ that is not isotopic mod
  $\t^1$ to a sub-disc of $S_0$, in contradiction to the preceding lemma.
  Thus by cut-and-paste along sub-discs of $S\setminus \t^1$, we can
  assume that additionally $S\cap \varphi(X\times I) =
  \varphi(X\times 1)$. 

  Let $\gamma\subset X$ be a circle so that $\varphi(\gamma\times 0)$ is
  contained in the boundary of a tetrahedron of $\t$.
  Since $S_0$ is 1--normal,  $\varphi(\gamma\times 0)$ bounds an open
  disc in  $S_0\setminus \t^2$. By the same argument as in the preceding
  paragraph, $\varphi(\gamma\times 1)$ bounds an open disc in $S\setminus
  \t^1$. 
  One easily verifies that these two discs together with
  $\varphi(\gamma\times I)$ bound a ball in $B^+(S)\cap  B^-(S_0)$
  disjoint to $\t^1$. 
  Hence $(B^+(S)\cap  B^-(S_0))\setminus U(\varphi(X\times I))$ is a
  disjoint union of balls in $M\setminus \t^1$, and this implies the
  existence of $J$.
\end{proof}

\section{Reduction of surfaces}
\label{sec:reduction}

 Let $M$ be a closed connected orientable 3--manifold with a
 triangulation $\t$. In this
 section, we show how to get isotopies 
 of embedded surfaces under which the number of
 intersections with $\t^1$ is monotonely non-increasing. 

\begin{definition}
  Let $S\subset M$ be a connected embedded surface that is transversal
  to $\t^2$ and splits $M$.
  Let $D$ be an upper (resp.\ lower) bond of $S$, set $D_1=U(D)\cap
  S$, and set $D_2= B^+(S)\cap \d U(D)$ (resp.\ $D_2= B^-(S)\cap \d
  U(D)$). 
  An \begriff{elementary reduction} along $D$ transforms $S$ to the
  surface $(S\setminus D_1)\cup D_2$. 
  \begriff{Upper} (resp.\ \begriff{lower}) \begriff{reductions} of $S$
  are the surfaces that are obtained from $S$ by
  a sequence of elementary reductions along upper (resp.\ lower)
  bonds. 
\end{definition}
If $S'$ is an upper or lower reduction of $S$, then $\|S'\|\le \|S\|$
with equality if and only if $S=S'$. Obviously $S$ is isotopic to
$S'$, such that $\|\cdot\|$ is monotonely non-increasing under the isotopy. 
If $\alpha\subset \t^1\setminus \t^0$ is an arc
with $\d\alpha\subset S'$, then also $\d\alpha\subset S$.
It is easy to see that if $S'$ has a lower compressing disc and is an
upper reduction of $S$, then also $S$ has a lower compressing
disc.

We will construct surfaces
with almost 1--normal upper or lower reductions. 
Let $N\subset M$ be a
3--dimensional sub--manifold, such that $\d N$ is pre-normal. 
Let $S\subset N$ be an embedded surface transversal to
$\t^2$ that splits $M$ and has no lower compressing disc.  

\begin{lemma}\label{lem:r+a1n}
  Suppose that there is a system $\Gamma\subset N\setminus \t^1$ of 
  arcs such that $S^\Gamma\subset N$ is connected, $\Gamma\subset
  B^-(S^\Gamma)$,  and $\d N\cap
  B^+(S^\Gamma)$ is 1--normal. 

  If, moreover, $\Gamma$ and an upper reduction $S'\subset N$ of
  $S^\Gamma$ are chosen so that $\|S'\|$ is minimal, then $S'$ is almost
  1--normal. 
\end{lemma}
\begin{proof}
  By hypothesis, $\Gamma\subset B^-(S^\Gamma)$, and $S$ 
  has no lower compressing discs. Thus by Lemma~\ref{lem:keinelcd}, 
  $S^\Gamma$ has no lower compressing discs. Therefore its upper reduction
  $S'$ has no lower compressing discs.

  Assume that $S'$ is not almost 1--normal. 
  Then $S'$ has a compressing disc $D'$ that is contained in a single
  tetrahedron $t$ (see~\cite{matv}), with string $\alpha'$ and base
  $\beta'$. Since $S'$ has no lower compressing discs, $D'$ is upper
  and does not contain proper compressing sub-discs. Thus $\alpha'\cap
  S'=\d\alpha'$, i.e., all components of $(D'\cap 
  S')\setminus \beta'$ are circles. 
  Since $\d N$ is pre-normal, $\d N\setminus \t^2$ is a disjoint union
  of discs. Therefore, since $D'$ is contained in a single tetrahedron,
  we can assume by isotopy of $D'$ mod $\t^2$ that $D'\cap \d 
  N$ consists of arcs. 
  We have $\alpha'\subset B^+(S')\subset B^+(S^\Gamma)$. It   follows $\d
  N\cap \alpha'=\emptyset$, since otherwise a sub-disc of $D'$ is a
  compressing disc for $\d N\cap 
  B^+(S^\Gamma)$, which is impossible as $\d N\cap B^+(S^\Gamma)$ is 1--normal by
  hypothesis. Thus $\d N\cap \alpha'=\emptyset$ and $D'\subset N$. 
\begin{figure}[ht!]
\cl{\relabelbox\small
\epsfbox{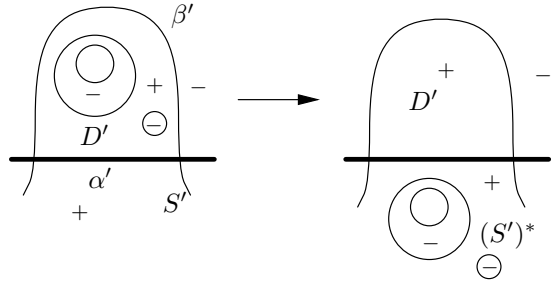}\let\ss\scriptsize
\relabel {D}{$D'$}
\relabel {D1}{$D'$}
\relabel {a}{$\alpha'$}
\relabel {b}{$\beta'$}
\relabela <-1pt,0pt> {S}{$S'$}
\relabela <0pt,1pt> {S1}{$(S')^*$}
\relabela <-1pt,0pt> {+}{\ss$+$}
\relabela <-1pt,0pt> {+1}{\ss$+$}
\relabela <-1pt,0pt> {+2}{\ss$+$}
\relabela <-1pt,0pt> {+3}{\ss$+$}
\relabela <-2pt,0pt> {-}{\ss$-$}
\relabela <-2pt,0pt> {-1}{\ss$-$}
\relabela <-2pt,0pt> {-2}{\ss$-$}
\relabela <-2pt,0pt> {-3}{\ss$-$}
\relabela <-2pt,0pt> {-4}{\ss$-$}
\relabela <-2pt,0pt> {-5}{\ss$-$}
\endrelabelbox}
\caption{How to produce a bond}
\label{fig:makebridge}
\end{figure}

  By an isotopy with support in $U(D')$ that is constant on
  $\beta'$, we move $(D'\cap S')\setminus \beta'$ to $U(D')\setminus
  t$, and obtain from $S'$ a surface $(S')^*\subset N$ that has $D'$
  as upper bond. This is shown in Figure~\ref{fig:makebridge}, where
  $B^+(S')$ is indicated by plus signs and $\t^1$ is bold. The
  isotopy moves $\Gamma$  
  to a system of arcs $\Gamma^*\subset N$ and moves $S^\Gamma$ to 
  $S^{\Gamma^*}$ with $\Gamma^*\subset B^-(S^{\Gamma^*})$. 
  Since $\alpha'\subset B^+(S')$, there is a
  homeomorphism $\varphi\co  B^-(S')\to B^-((S')^*)$ that is constant on
  $\t^1$ with $\varphi(B^-(S^\Gamma))=B^-(S^{\Gamma^*})$. 
  One obtains $S'$
  by a sequence of elementary reductions 
  along bonds of $S^\Gamma$ that are contained in $B^-(S')$. These bonds
  are carried by $\varphi$ to bonds of $S^{\Gamma^*}$. Thus $(S')^*$ is an
  upper reduction of $S^{\Gamma^*}$. Since $(S')^*$ admits an elementary
  reduction along its upper bond $D'$, we obtain a contradiction to
  the minimality of $\|S'\|$. Thus $S'$ is almost 1--normal.
\end{proof}

\begin{lemma}\label{lem:r+a1n*}
  Let $\Gamma$ and $S'$ be as in the previous lemma, and let $G_1, G_2$
  be two connected 
  components of $(S')^\times$ that both split $M$.
  Then there is no
  arc in $(\t^1\setminus\t^0)\cap B^+(S')\cap N$ joining 
  $G_1$ with $G_2$.  
\end{lemma}
\begin{proof}
  By the previous lemma, $S'$ is almost 1--normal. Recall that one
  obtains $(S')^\times$ up to isotopy mod $\t^1$ by splitting $S'$
  along splitting discs that do not meet $\t^2$. Assume that there
  is an arc $\alpha\subset(\t^1\setminus\t^0)\cap 
  B^+(S')\cap N$ joining $G_1$ with $G_2$. Let $Y$ be the
  component of $M\setminus (G_1\cup G_2)$ that contains
  $\alpha$.

  By hypothesis, $S^\Gamma$ is connected. Thus $S'$ is connected, and there is
  an arc $\beta\subset S'$ with 
  $\d\beta=\d\alpha$. Since $G_1,G_2$ split $M$, the set $Y$ is the
  only component of $M\setminus (G_1\cup G_2)$ with boundary $G_1\cup
  G_2$. Thus there is  a component $\beta'$ of
  $\beta\cap Y$ connecting $G_1$ with $G_2$. There is a splitting disc
  $D\subset Y$ of $S'$ contained in a single tetrahedron with
  $\beta'\cap D\not= \emptyset$. By choosing $D$ innermost, we assume
  that $\beta\cap D$ is a single point in $\d D$. Since $\d N$ is
  pre-normal and $D$ is contained in a single tetrahedron,
  we can assume by isotopy of $D$ mod $\t^2$ that $D\cap \d
  N=\emptyset$, thus $D\subset N$.

  Choose a disc $D'\subset U(\alpha\cup\beta)\cap B^+(S')$ so that
  $D'\cap \t^1=\alpha$ and $D'\cap S'=\beta\setminus U(\d D)$. Then 
  $D'\cap\d N=\emptyset$, since $U(\alpha\cup\beta)\cap \d
  N=\emptyset$. We split $S'$ along $D$, pull the two components of
  $(S'\cap \d U(D))\setminus D$ along  $(\d
  D')\setminus (\alpha\cup\beta)$, and reglue. We obtain a surface
  $(S')^*$ with $D'$ as an upper bond. 

  Since a small collar of $\d D$ in $D$ is in $B^-(S')$, there is a
  homeomorphism $\varphi\co  B^-(S')\to B^-((S')^*)$ that is constant on
  $\t^1$. Set $\Gamma^*=\varphi(\Gamma)$.
  Then $\varphi(S^\Gamma) = S^{\Gamma^*}$ with
  $\Gamma^*\subset B^-(S^{\Gamma^*})$. As in the proof of the previous lemma,
  $(S')^*$ is an upper reduction of $S^{\Gamma^*}$, and $(S')^*$ admits an
  elementary reduction along $D'$.
  This contradiction  to the minimality of $\|S'\|$ yields the lemma. 
\end{proof}


\section{Proof of Theorem~\ref{prop:}}
\label{sec:theproof}

Let $\t$ be a triangulation of $S^3$ with a  vertex
$x_0\in\t^0$. Let $\Sigma\subset S^3$ be a maximal system of
disjoint 1--normal spheres with $\| \Sigma\|< 2^{185 \,t(\t)^2}$, as
given by Construction~\ref{con:Sigma}. 
Construction~\ref{con:tildeSigma} extends $\Sigma$ to a system
$\tilde\Sigma\subset S^3$ of disjoint 2--normal spheres that are
pairwise non-isotopic mod $\t^2$,  
such that  
\begin{enumerate}
\item any component of $\tilde\Sigma$ has at most one octagon,
\item any component of $S^3\setminus \tilde\Sigma$ has at most one
  boundary component that is not 1--normal,
\item if the boundary of a component $N$ of $S^3\setminus \tilde
  \Sigma$ is 1--normal, then $N$ does not contain 2--normal spheres with
  exactly one octagon, and
\item $\|\tilde\Sigma\| < 2^{190 \,t(\t)^2}$.
\end{enumerate}

Let $N$ be a component of $S^3\setminus \tilde\Sigma$ that is not a
regular neighbourhood of a vertex of $\t$. 
Let $S_0$ be the component of $\d N$ with $N\subset
B^-(S_0)$, and let $S_1,\dots, S_k$ be the other components of $\d
N$. Since $\Sigma$ is maximal, any almost 1-normal sphere in $N$ is a 
tube sum of copies of $S_0,S_1,\dots, S_k$. 
\begin{lemma}\label{lem:Nohnev}
  $N\cap \t^0=\emptyset$.
\end{lemma}
\begin{proof}
  If $x\in N\cap\t^0$, then the sphere $\d U(x)\subset N$ is
  1--normal. It is not isotopic mod $\t^1$ to a component of $\d N$,
  since $N\not= U(x)$. This contradicts the maximality of
  $\Sigma$. 
\end{proof}

\begin{lemma}\label{lem:gibtFinN}
  If $\d N$ is 1--normal, then there is an arc in $\t^1\cap \overline
  N$ that connects two different components of $\d N\setminus S_0$.
\end{lemma}
\begin{proof}
  Let $\d N=S_0\cup S_1\cup\dots\cup S_k$ be 1--normal. We first consider the
  case where there is an almost 1--normal sphere $S\in S_1\circ
  \dots\circ S_k$ in $\overline N$ that has a compressing disc $D$, with
  string $\alpha$ and base $\beta$. 
  We choose $D$ innermost, so that 
  $\alpha\cap S=\d\alpha$. In particular, $\alpha\cap \d N=\d
  \alpha$. Assume that $\alpha\not\subset \overline N$. Since $\d
  D\setminus\alpha\subset \overline N$, there 
  is an arc $\beta'\subset D\cap \d N$  that connects the endpoints of
  $\alpha$. The sub-disc $D'\subset D$ bounded by
  $\alpha\cup\beta'$ is a compressing disc for the 1--normal surface $\d
  N$, in contradiction to Lemma~\ref{lem:1ohnedisc}. By consequence,
  $\alpha\subset \overline N$. 
  Assume that $\d\alpha$ is contained in a single component of $\d
  N\setminus S_0$, say, in $S_1$. By Lemma~\ref{lem:1ohnedisc}, $D$ is
  not a compressing disc for $S_1$, hence $\beta\not\subset
  S_1$.
  Thus there is a closed line in $S_1\setminus \beta$
  that separates $\d\alpha$ on $S_1$, but not on $S$.  This is
  impossible as $S$ is a sphere. We conclude that if $S$ has a
  compressing disc, then there is an arc $\alpha\subset \t^1\cap N$
  that connects different components of $\d N\setminus S_0$.

  It remains to consider the case where no sphere in $S_1\circ
  \dots\circ S_k$ contained in $\overline N$ has a compressing disc.
  We will show the existence of an almost 2--normal sphere in $N$ with
  exactly one octagon, using the technique of thin position. 
  This contradicts property (3) of 
  $\tilde \Sigma$ (see the begin of this section), and
  therefore finishes the proof of the lemma. 
  Let $J\co S^2\times I\nach B^-(S_0)$ be a $\t^1$--Morse embedding, such
  that  
  \begin{enumerate}
  \item $J(S^2\times 0)=S_0$,
  \item $J(S^2\times \frac 12)\in S_1\circ \dots\circ S_k$  
    (or $\|J(S^2\times\frac 12)\|=0$, in the case $\d N=S_0$), 
  \item $B^-\left(J(S^2\times 1)\right)\cap \t^1=\emptyset$, and
  \item $\kappa(J)$ is minimal. 
  \end{enumerate}
  Define $S=J(S^2\times\frac 12)$. 
  Assume that for some $\xi\in I$ there is a pair $D_1,D_2\subset M$
  of nested or independent upper and lower compressing discs for
  $J_\xi=J(S^2\times \xi)$. 
  We show that we can assume $D_1,D_2\subset B^-(S_0)$. 
  Since $S_0$ is
  1--normal, it has no compressing 
  discs by Lemma~\ref{lem:1ohnedisc}. Thus $(D_1\cup D_2)\cap S_0$
  consists of circles. Any such 
  circle bounds a disc in $S_0\setminus \t^1$ by
  Lemma~\ref{lem:discparallel}. By cut-and-paste of $D_1\cup D_2$, we
  obtain $D_1,D_2\subset B^-(S_0)$, as claimed.  
  Now, one obtains from $J$ an embedding $J'\co S^2\times 
  I\nach B^-(S_0)$ with $\kappa(J')<\kappa(J)$ by isotopy along
  $D_1\cup D_2$, see~\cite{matv}, \cite{thompson}. The embedding
  $J'$ meets conditions (1) and (3) in the definition of
  $J$. 
  Since $S\in S_1\circ\dots \circ S_k$ has no compressing discs by assumption,
  $S\cap D_i$ consists of circles. Thus $S$ is split equivalent to 
  $J'(S^2\times \frac 12)$. 
  So $J'$ meets also condition (2), $J'(S^2\times \frac 12)\in
  S_1\circ \dots\circ S_k$, in  
  contradiction to the choice of $J$. This disproves the existence of
  $D_1,D_2$. In conclusion, if $J_\xi$ has upper and lower bonds,
  then it is impermeable.

  Let $\xi_{max}$ be the greatest critical parameter of
  $J$ with respect to $\t^1$ in the interval $\left]0,\frac
    12\right[$. 
  We have $N\cap \t^0=\emptyset$ by Lemma~\ref{lem:Nohnev}.
  Hence the critical point corresponding to $\xi_{max}$ is a point of
  tangency of  $J_{\xi_{max}}$ to some edge of 
  $\t$. By assumption, $S$ has no upper bonds, thus
  $\|S\|<\|J_{\xi_{max}-\epsilon}\|$ for sufficiently
  small $\epsilon>0$.
  Let $\xi_{min}\in I$ be the smallest critical parameter of $J$ with
  respect to $\t^1$. By Lemma~\ref{lem:1ohnedisc}, $S_0$ has no 
  bonds, thus $\|S_0\|<\|J_{\xi_{min}+\epsilon}\|$.
  Therefore there are consecutive critical parameters
  $\xi_1,\xi_2\in\left]0,\frac 12\right[$ such that 
  $$ \|J_{\xi_1-\epsilon}\| <
     \|J_{\xi_{1}+\epsilon}\| >
     \|J_{\xi_{2}+\epsilon}\|.$$
  Thus $J_{\xi_{1}+\epsilon}$ has both upper and lower
  bonds, and is therefore impermeable by the preceding
  paragraph. One component of $J_{\xi_{1}+\epsilon}^\times$ is a
  2--normal sphere in $N$ with exactly one octagon, by
  Proposition~\ref{prop:impermeabel}. The existence of that 2--normal
  sphere is a contradiction to the properties of $\tilde\Sigma$, which
  proves the lemma.  
\end{proof}

We show that some tube sum $S\in S_1\circ\dots\circ S_k$ is isotopic
to $S_0$ such that $\|\cdot\|$ is monotone under the isotopy. We
consider three cases. In the first case, let $\d N$ be 1--normal. 
\begin{lemma}\label{lem:r+S}
  If $\d N$ is 1--normal, then there is a sphere $S\in S_1\circ
  \dots\circ S_k$ in $N$ with an upper reduction $S'\subset N$ so that
  there is a $\t^1$--Morse embedding $J\co S^2\times I \to S^3$ with 
  $J(S^2\times I) = B^+(S') \cap B^-(S_0)$ and $c(J,\t^1)= 0$.  
\end{lemma}
\begin{proof}
  By Lemma~\ref{lem:gibtFinN}, there is an arc $\alpha\subset
  \t^1\cap N$ that connects two components of $\d N\setminus S_0$,
  say, $S_1$ with $S_2$. By Lemma~\ref{lem:Nohnev}, $\alpha$ is
  contained in an edge of $\t$.
  By Lemma~\ref{lem:1ohnedisc},
  the 1--normal surfaces $S_1,\dots,S_k$ have no lower compressing
  discs. 
  Let $\Gamma\subset N$ be a system of $k-1$ arcs, such that the tube 
  sum $S$ of $S_1,\dots,S_k$ along $\Gamma$ is a sphere and an
  upper reduction $S'\subset N$ of $S$ minimizes $\|S'\|$. We have
  $\|S'\|<\|S\|$, since it is possible to choose $\Gamma$ so that
  $S$ has an upper bond with string $\alpha$.
  Since $\Gamma\subset B^-(S)$ and by Lemma~\ref{lem:r+a1n}, $S'$ is
  almost 1--normal. 

  By the maximality of $\Sigma$, it follows $S'\in
  n_0S_0\circ \dots\circ  n_kS_k$ with non-negative integers
  $n_0,n_1,\dots,n_k$. Moreover, $n_i\le 2$ for $i=0,\dots,k$ by
  Lemma~\ref{lem:r+a1n*}. Since $S$ separates $S_0$
  from $S_1,\dots,S_k$, so does $S'$. Thus any path connecting $S_0$
  with $S_j$ for some $j\in\{1,\dots,k\}$ intersects $S'$ in an odd
  number of points. So alternatively $n_0\in \{0,2\}$ and $n_i=1$ for
  all $i\in\{1,\dots,k\}$, or $n_0=1$ and $n_i\in\{0,2\}$ for all
  $i\in\{1,\dots,k\}$. Since $\|S'\|< \|S^*\|$, it follows
  $n_0=1$ and $n_i=0$ for $i\in\{1,\dots,k\}$, thus $(S')^\times=S_0$. 
  The existence of a $\t^1$--Morse
  embedding $J$ with the claimed properties follows then by
  Corollary~\ref{cor:gibtJ}. 
\end{proof}

The second case is that $S_0$ is 1--normal, and exactly one of
$S_1,\dots,S_k$ contains exactly one octagon, say, $S_1$. The octagon
gives rise to an upper bond $D$ of $S_1$ contained in a single
tetrahedron. Since $\d N\setminus S_1$ is 1--normal, $D\subset N$. Thus
an elementary reduction of $S_1$ 
along $D$ transforms $S_1$ to a sphere $F\subset N$. Since $S_1$ is
impermeable by Proposition~\ref{prop:2n-imperm}, $F$ has no lower
compressing disc (such a disc would give rise to a lower
compressing disc for $S_1$ that is independent from $D$).
\begin{lemma}\label{lem:r+F}
  If $\d N\setminus S_0$ is not 1--normal, then there is a sphere $S\in
  S_1\circ \dots\circ S_k$ in $N$ with an upper reduction $S'\subset N$
  so that  there is a $\t^1$--Morse embedding $J\co S^2\times I \to S^3$ with 
  $J(S^2\times I) = B^+(S') \cap B^-(S_0)$ and $c(J,\t^1)= 0$.
\end{lemma}
\begin{proof}
  We apply the Lemma~\ref{lem:r+a1n} to $F,S_2,\dots,S_k$, and 
  together with the elementary reduction along $D$ we obtain a sphere
  $S\in S_1\circ S_2\circ\dots\circ S_k$ with an almost 1--normal
  upper reduction $S'\subset N$.  
  One concludes $(S')^\times = S_0$ and the existence of $J$ as in the
  proof of the   previous lemma. 
\end{proof}

We come to the third and last case, namely $S_0$ has exactly one
octagon and $\d N\setminus S_0$ is 1--normal. The octagon gives rise to 
a lower bond $D$ of $S_0$, that is contained in $N$ since $\d
N\setminus S_0$ is 1--normal. Thus an elementary reduction of $S_0$
along $D$ yields a sphere $F\subset N$. Since $S_0$ is impermeable by
Proposition~\ref{prop:2n-imperm}, $F$ has no upper compressing disc,
similar to the previous case. 
\begin{lemma}\label{lem:r-F}
  If $S_0$ is not 1--normal, then there is a lower reduction $S'\in
  S_1\circ \dots\circ S_k$ of $S_0$, with $S'\subset N$. 
\end{lemma}
\begin{proof}
  We apply Lemma~\ref{lem:r+a1n} with $\Gamma=\emptyset$ to
  \emph{lower} reductions of $F$, which is possible by symmetry.
  Thus, together with the elementary reduction along
  $D$, there is a lower reduction $S'\in n_0S_0\circ \dots \circ
  n_kS_k$ of $S_0$, and $n_0,\dots,n_k\le 2$ by
  Lemma~\ref{lem:r+a1n*}. Since $S'\subset B^-(F)$ and 
  $S_0\subset B^+(F)$, it follows $n_0=0$. 
  Since $S'$ separates $\d N\cap B^+(F)$ from $\d N\cap B^-(F)$, it
  follows $n_1,\dots,n_k$ odd, thus $n_1=\dots = n_k=1$.  
\end{proof}

We are now ready to construct the $\t^1$--Morse embedding $H\co S^2\times
I\to S^3$ with $c(H,\t^1)$ bounded in terms of $t(\t)$, thus to finish
the proof of Theorems~\ref{thm:main} and~\ref{prop:}. 
Let $x_0\in \t^0$ be the vertex involved in the definition of
$B^+(\cdot)$. We construct $H$ inductively as follows.

Choose $\xi_1\in \left]0,1\right[$ 
and choose $H|[0,\xi_1]$ so that $H_0\cap \t^2=\emptyset$,
$H_{\xi_1} = \d U(x_0)\subset \tilde\Sigma$, and $x_0$ is the only
critical point of $H|[0,\xi_1]$.  

For $i\ge 1$, let $H|[0,\xi_i]$ be already constructed. 
Our induction hypothesis is that 
$H_{\xi_i}\in S_0\circ S^*$ for some component $S_0$ of
$\tilde\Sigma$, and moreover for any choice of $S_0$ we have $H_{\xi_i}
\subset B^+(S_0)$. 
Choose $\xi_{i+1}\in \left]\xi_i, 1\right[$. 

Assume that $S_0$ is not of the form $S_0=\d U(x)$ for a vertex $x\in
\t^0\setminus \{x_0\}$. 
Then, let $N_i$ be the component of $S^3\setminus \tilde\Sigma$ with
$N_i\subset B^-(S_0)$ and $\d N_i= S_0\cup S_1\cup\dots\cup S_k$ for
$S_1,\dots,S_k\subset \tilde\Sigma$. 
If $S_0$ is 1--normal, then let $S\in S_1\circ\dots\circ S_k$, $S'$ and
$J$ be as in Lemmas~\ref{lem:r+S} and~\ref{lem:r+F}.
Then, we extend $H|[0,\xi_i]$ to $H|[0,\xi_{i+1}]$ induced by the
embedding $J$, relating $S_0$ with $S'$, and 
by the \emph{inverses} of the elementary upper reductions, relating $S'$
with $S$. 
If $S_0$ is not 1--normal, then let $S\in S_1\circ\dots\circ S_k$ be as
in Lemma~\ref{lem:r-F}. 
We extend $H|[0,\xi_i]$ to $H|[0,\xi_{i+1}]$ along the elementary lower
reductions, relating $S_0$ with $S$.
In either case, $H_{\xi_{i+1}} \in S_1\circ\dots \circ S_k\circ S^*$.
The critical points of $H|[\xi_i,\xi_{i+1}]$ are contained in $N_i$, given
by elementary reductions.  
Thus the number of these critical points is $\le
\frac 12 \max \{\|S_0\|, \|S\|\}\le \frac 12\, \|\tilde \Sigma\| <
2^{190 \,t(\t)^2}$, by Construction~\ref{con:tildeSigma}.  
Since 
$H_{\xi_{i+1}} \subset B^+(S_m)$ for any $m=1,\dots, k$, we can proceed
with our induction. 

After at most $\#(\tilde \Sigma)$ steps, we have $H_{\xi_i}^\times = \d 
U(\t^0\setminus \{x_0\})$. Then, choose $H|[\xi_{i},1]$ so that $H_1\cap
\t^2=\emptyset$ and the set of its critical points is $\t^0\setminus
\{x_0\}$. 
By Proposition~\ref{lem:kneserhaken} holds $\#(\tilde \Sigma) \le 10 \,t(\t)$.
Thus finally $$c(H,\t^1) < \#(\t^0) + 10 \,t(\t)\cdot
2^{190 \,t(\t)^2} < 2^{196 \,t(\t)^2}.\qed$$


\section{Proof of Propositions~\ref{prop:impermeabel}
  and~\ref{prop:2n-imperm}}
\label{sec:imperm-a2}

Let $M$ be a closed connected 3--manifold with a triangulation
$\t$. We 
prove Proposition~\ref{prop:impermeabel}, that states that any
impermeable surface in $M$ is isotopic mod $\t^1$ to an almost
2--normal surface with exactly one octagon. The proof consists of the
following three lemmas. 

\begin{lemma}
  Any impermeable surface in $M$ is almost 2--normal, up to isotopy mod
  $\t^1$.
\end{lemma}
\begin{proof}
  We give here just an outline. A complete proof can be found
  in~\cite{matv}. Let $S\subset M$ be an impermeable surface. By
  definition, it has upper and lower bonds with strings
  $\alpha_1,\alpha_2$. By isotopies mod $\t^1$, one obtains
  from $S$ two surfaces $S_1,S_2\subset M$, such that $S_i$ 
  has a return $\beta_i\subset \t^2$ with $\d\beta_i=\d\alpha_i$,
  for $i\in\{1,2\}$. A surface that has both upper and lower returns
  admits an independent pair of upper and lower compressing discs,
  thus is not impermeable. 
  By consequence, under the isotopy mod $\t^1$ that relates
  $S_1$ and $S_2$ occurs a surface $S'$ that has no
  returns at all, thus is almost $k$--normal for some natural number
  $k$. 

  If there is a boundary component $\zeta$ of a component of
  $S'\setminus \t^2$ and an edge $e$ of $\t$ with $\#(\zeta\cap e)>2$,
  then there is an independent pair of upper and lower compressing
  discs. Thus $k=2$.
\end{proof}

\begin{lemma}\label{lem:2octagon}
  Let $S\subset M$ be an almost 2--normal impermeable surface. Then $S$
  contains at most one octagon. 
\end{lemma}
\begin{proof}
  Two octagons in different tetrahedra of $\t$ give rise to a pair of
  independent upper and lower compressing discs for $S$. Two octagons
  in one tetrahedron of $\t$ give rise to a pair of nested upper
  and lower compressing discs for $S$. Both is a contradiction to the
  impermeability of $S$. 
\end{proof}

\begin{lemma}\label{lem:1nichtimpermeabel}
  Let $S\subset M$ be an almost 2--normal impermeable surface. Then $S$
  contains at least one octagon. 
\end{lemma}
\begin{proof}
  By hypothesis, $S$ has both upper and lower bonds. Assume that $S$
  does not contain octagons, i.e., it
  is almost 1--normal. We will obtain a contradiction to the
  impermeability of $S$ by showing that $S$ has a pair of independent or
  nested compressing discs.  

  According to Lemma~\ref{lem:1ohnedisc}, we can assume that $S$ has a
  compressing disc $D_1$ with string $\alpha_1$ that is contained in a
  single closed tetrahedron $t_1$. Choose $D_1$ innermost,
  i.e., $\alpha_1\cap S=\d\alpha_1$. Without assumption, let $D_1$ be 
  \emph{upper}. Since $S$ has no octagon by assumption,
  $\alpha_1$ connects two different components $\zeta_1,\eta_1$ of
  $S\cap \d t_1$. Let $D$ be a lower bond of $S$. Choose $S$, $D_1$
  and $D$ so that, in addition, $\#(D\cap\t^2)$ is minimal.  

  Let $C$ be the closure of an innermost component of 
  $D\setminus\t^2$, which is a disc. There is a closed tetrahedron 
  $t_2$ of $\t$ and a closed 2--simplex $\sigma_2\subset\d t_2$ of $\t$
  such that $\d C\cap \d t_2$ is a single component $\gamma\subset 
  \sigma_2$. We have to consider three cases.
  \begin{enumerate}
  \item Let $\gamma$ be a circle, thus $\d C=\gamma$. There is a
    disc $D'\subset\sigma_2$ with $\d D'=\gamma$ and a ball
    $B\subset t_2$ with $\d B=C\cup D'$. We move $S\cup D$ away from
    $B$ by an isotopy mod $\t^1$ with support in
    $U(B)$, and obtain a surface $S^*$ with a lower bond $D^*$. As
    $D$ is a bond, $S\cap D'$ consists of circles. Therefore the
    normal arcs of $S\cap \t^2$ are not changed under the isotopy, and
    the isotopy does not introduce returns, thus $S^*$ is almost
    1--normal. Since $\xi_1\cap D'=\eta_1\cap D'=\emptyset$ and $C\cap
    S=\emptyset$, it follows  $B\cap \d D_1=\emptyset$. Thus $D_1$ is
    an upper compressing disc for $S^*$, and
    $\#(D^*\cap \t^2)<\#(D\cap \t^2)$ in 
    contradiction to our choice.    
  \item Let $\gamma$ be an arc with endpoints in a single component
    $c$ of $S\cap\sigma$. By an isotopy mod $\t^1$
    with support in $U(C)$ that moves $C$ into
    $U(C)\setminus t_2$, we obtain from $S$ and $D$ a surface $S^*$
    with a lower bond $D^*$. Since $D$ is a bond, the isotopy does
    not introduce returns, thus $S^*$ is almost 1--normal. One
    component of $S^*\cap t_1$ is isotopic mod $\t^2$ to the
    component of $S\cap t_1$ that contains $\d D_1\cap S$. Thus up to
    isotopy mod $\t^2$,  $D_1$ is an upper compressing disc for
    $S^*$, and $\#(D^*\cap \t^2)<\#(D\cap \t^2)$ in 
    contradiction to our choice.    
  \item  Let $\gamma$ be an arc with endpoints in two different
    components $c_1,c_2$ of $S\cap\sigma$.  
    Assume that, say, $c_1$ is a circle. By an isotopy mod $\t^1$
    with support in $U(C)$ that moves $C$ into
    $U(C)\setminus t_2$, we obtain from $S$ and $D$ a surface $S^*$
    with a lower bond 
    $D^*$. Since $D$ is a bond, the isotopy does not introduce
    returns, thus $S^*$ is almost 1--normal. There is a disc
    $D'\subset \sigma$ with $\d D'=c_1$. Let $K$ be the component of
    $S\cap t_1$ that contains $\d D_1\cap S$. One component of
    $S^*\cap t_1$ is isotopic mod $\t^2$ either to $K$ or, if
    $\d D'\cap\d K\not=\emptyset$, to $K\cup D'$. In either case,
    $D_1$ is an upper  compressing disc for $S^*$, up to isotopy
    mod $\t^2$. But $\#(D^*\cap \t^2)<\#(D\cap \t^2)$ in  
    contradiction to our choice. Thus, $c_1$ and $c_2$ are normal
    arcs.
  \end{enumerate}
  Since $S$ is almost 1--normal, $c_1$, $c_2$ are contained in
  different components $\zeta_2,\eta_2$ of $S\cap \d t_2$.
  Since $D$ is a lower bond, $\d(C\cap D_1)\subset \d C\cap
  S$. There is a sub-arc $\alpha_2$ of an edge of $t_2$ 
  and a disc $D'\subset \sigma$ with $\d D'\subset \alpha_2\cup
  \gamma\cup \zeta_2\cup\eta_2$ and $\alpha_2\cap S=\d\alpha_2$. The
  disc $D_2=C\cup D'\subset t_2$ is a lower compressing disc for $S$
  with string $\alpha_2$, and $\d(D_1\cap D_2)\subset
  \d D_2\cap S$.
  At least one component of $\d t_1\setminus(\zeta_1 \cup \eta_1)$ is a 
  disc that is disjoint to $D_2$. Let $D_m$ be the closure of a
  copy of such a disc in the interior of $t_1$, with $\d D_m\subset
  S$. By construction, $D_1\cap D_m=\d
  D_1\cap \d D_m$ is a single point and $D_2\cap D_m=\emptyset$. Thus
  by Lemma~\ref{lem:makeindep}, $S$ has a pair of independent or nested
  upper  and lower compressing discs and is therefore not impermeable.
\end{proof}

\paragraph*{Proof of Proposition~\ref{prop:2n-imperm}}
Let $S\subset M$ be a connected 2--normal
surface  that splits $M$, and assume that exactly one
component $O$ of $S\setminus \t^2$ is an octagon. The octagon gives
rise to upper and lower bonds of $S$. 

Let $D_1,D_2$ be any upper and lower compressing discs for $S$. We
have to show that $D_1$ and $D_2$ are neither impermeable nor
nested. It suffices to show that $\d D_1\cap \d
D_2\not\subset \t^1$.  
To obtain a contradiction, assume that $\d D_1\cap \d D_2\subset
\t^1$. Choose $D_1,D_2$ so that $\#(\d D_1\setminus \t^2) +
\#(\d D_2\setminus \t^2)$ is minimal. 

Let $t$ be a tetrahedron of $\t$ with a closed 2--simplex
$\sigma\subset \d t$, and let $\beta$ be a component of $\d D_1\cap t$
(resp.\ $\d D_2\cap t$) such that $\d \beta$ is
contained in a single component of $S\cap \sigma$. Since $S$ is
2--normal, there is a disc $D\subset S\cap t$ and an arc $\gamma\subset
S\cap \sigma$ with $\d D=\beta\cup \gamma$. By choosing $\beta$
innermost in $D$, we can assume that $D\cap (\d D_1\cup\d D_2) =
\beta$. An isotopy of $(D_1,\d D_1)$ (resp.\ $(D_2,\d D_2)$) in
$(M,S)$ with support in $U(D)$ that moves $\beta$ to $U(D)\setminus t$ 
reduces $\#(\d D_1\setminus \t^2)$ (resp.\  $\#(\d D_2\setminus
\t^2)$), leaving $\d D_1\cap \d D_2$ unchanged. This is a 
contradiction to the minimality of $D_1,D_2$.

For $i=1,2$, there are arcs
$\beta_i\subset \d D_i\setminus \t^1$ and $\gamma_i\subset 
D_i\cap \t^2$ such that $\beta_i\cup\gamma_i$ bounds a component of
$D_i\setminus \t^2$, by an innermost arc argument. Let $t_i$
be the tetrahedron of $\t$ that contains $\beta_i$, and let
$\sigma_i\subset \d t_i$ be the close 2--simplex that contains
$\gamma_i$. 
We have seen above that $\d \beta_i$ is not contained in a single
component of $S\cap \sigma_i$. Since $S$ is 2--normal, i.e., has no 
tubes, it follows that $\beta_i\subset O$. Since collars of $\beta_1$
in $D_1$ and of $\beta_2$ in $D_2$ are in different components of
$t\setminus O$, it follows $\beta_1\cap
\beta_2\not=\emptyset$. Thus $\d D_1\cap \d D_2\not\subset 
\t^1$, which yields Proposition~\ref{prop:2n-imperm}.\qed


\end{document}

%% file: gtoutput.tex
\def\ifplaintex{\expandafter\ifx\csname documentclass\endcsname\relax}

\ifplaintex 
\else
\fi

\expandafter\ifx\csname beginpicture\endcsname\relax
\expandafter\ifx\csname documentclass\endcsname\relax
\input pictex \else
\input prepictex \input pictex \input postpictex \fi\fi

\def\gt{{\mathsurround=0pt\it $\cal G\mskip-2mu$eometry \&\ 
$\cal T\!\!$opology}}        

\def\gtp{{\mathsurround=0pt\it $\cal G\mskip-2mu$eometry \&\ 
$\cal T\!\!$opology $\cal P\!$ublications}}  

\def\lognumber#1{\def\thelognumber{#1}}
\def\volumenumber#1{\def\thevolumenumber{#1}}
\def\papernumber#1{\def\thepapernumber{#1}}
\def\volumeyear#1{\def\thevolumeyear{#1}}

\def\pagenumbers#1#2{\def\startpage{#1}\def\finishpage{#2}}
\def\published#1{\def\publishdate{#1}}
\def\proposed#1{\def\theproposer{#1}}
\def\seconded#1{\def\theseconders{#1}}
\def\received#1{\def\receiveddate{#1}}

\def\accepted#1{\def\accepteddate{#1}}

\def\asciiaddress#1{\def\theasciiaddress{#1}}

\long\def\asciiabstract#1{\long\def\theasciiabstract{#1}}
\def\asciikeywords#1{\def\theasciikeywords{#1}}

\let\\\par\let\thelognumber\relax
\let\thevolumenumber\relax\let\thepapernumber\relax
\let\thevolumeyear\relax\let\thesamplenumber\relax\let\startpage\relax
\let\finishpage\relax\let\publishdate\relax\let\receiveddate\relax
\let\reviseddate\relax\let\accepteddate\relax\let\theasciititle\relax
\let\theasciiauthors\relax\let\theasciiaddress\relax
\let\theasciiabstract\relax\let\theasciikeywords\relax
\let\theasciiemail\relax\let\theshortauthors\relax\let\theshorttitle\relax

\long\def\maketitlep{   

\count0=\startpage

\gt\hfill      
\beginpicture
\setcoordinatesystem units <0.33truein, 0.33truein> point at 2.2 0.9
\setplotsymbol ({$\cal G$})
\plotsymbolspacing=9truept
\circulararc 315 degrees from 0 1 center at 0 0
\setplotsymbol ({$\cal T$})
\circulararc 315 degrees from 1 -1 center at 1 0
\endpicture
\break
{\small\ifx\thesamplenumber\relax 
Volume \else Sample
\fi\thevolumenumber\ (\thevolumeyear)
\startpage--\finishpage\nl
Published: \publishdate}
\vglue 0.5truein plus 0.4fil minus 0.1truein

{\parskip=0pt\leftskip 0pt plus 1fil\def\\{\par\smallskip}{\ifplaintex\large
\else\Large\fi\bf\thetitle}\par\medskip}   

\vglue 0pt plus 0.1fil 

{\parskip=0pt\leftskip 0pt plus 1fil\def\\{\par}{\sc\theauthors}
\par\medskip}

\vglue 0pt plus 0.1fil 

{\small\parskip=0pt\let\newline\\
{\leftskip 0pt plus 1fil\def\\{\par}{\sl\theaddress}\par}
\expandafter\ifx\theemail\relax    
\relax\else\vglue 5pt plus 0.02fil minus 2pt\def\\{\stdspace{\rm 
and}\stdspace} 
\cl{Email:\stdspace\tt\theemail}\fi
\ifx\theurl\relax                  
\relax\else\vglue 5pt plus 0.02fil minus 2pt\def\\{\stdspace{\rm 
and}\stdspace}
\cl{URL:\stdspace\tt\theurl}\fi\par}

\vglue 7pt plus 0.3fil minus 3pt

{\bf Abstract}
\vglue 5pt plus 0.1fil minus 2pt

\theabstract

\vglue 7pt plus 0.3fil minus 3pt

{\bf AMS Classification numbers}\quad Primary:\quad \theprimaryclass

Secondary:\quad \thesecondaryclass

\vglue 5pt plus 0.3fil minus 2pt

{\bf Keywords}\quad \thekeywords

\vglue 10pt plus 0.5fil minus 5pt

{\small  Proposed: \theproposer\hfill Received: \receiveddate\nl
Seconded: \theseconders\hfill 
\ifx\reviseddate\relax                         
Accepted: \accepteddate                        
\else
Revised: \reviseddate                          
\fi}
\eject
}       

\let\maketitlepage\maketitlep
\let\maketitle\maketitlepage

\font\phead=cmsl9 scaled 950
\font=cmsl9 scaled 1050
\font\pnum=cmbx10 scaled 913
\font\lnum=cmbx10 
\font\pfoot=cmsl9 scaled 950
\font=cmsl9 scaled 1050
\ifplaintex
\headline{\vbox to 0pt{\vskip -4.5mm\line{\small\phead\ifnum
\count0=\startpage ISSN 1364-0380 (on line)
1465-3060 (printed) \hfill {\pnum\folio}\else\ifodd\count0\def\\{ }%
\ifx\theshorttitle\relax\thetitle\else\theshorttitle\fi\hfill{\pnum\folio}
\else\def\\{ and }{\pnum\folio}\hfill\ifx\theshortauthors\relax\theauthors
\else\theshortauthors\fi\fi\fi}\vss}}
\footline{\vbox to 0pt{\vglue 0mm\line{\small\pfoot\ifnum\count0=\startpage
\copyright\ \gtp\hfill\else
\gt, Volume \thevolumenumber\ (\thevolumeyear)\hfill\fi}\vss
}}
\else
\makeatletter
\def\@oddhead{{\small\lhead\ifnum\count0=\startpage ISSN 1364-0380 (on line)
1465-3060 (printed) \hfill {\lnum\number\count0}\else\ifodd\count0
\def\\{ }\ifx\theshorttitle\relax \thetitle \else\theshorttitle\fi\hfill
{\lnum\number\count0}\else\def\\{ and }{\lnum\number\count0}
\hfill\ifx\theshortauthors\relax 
\theauthors\else\theshortauthors\fi\fi\fi}}\def\@evenhead{\@oddhead}
\def\@oddfoot{\small\lfoot\ifnum\count0=\startpage\copyright\ \gtp\hfill\else
\gt, Volume \thevolumenumber\ (\thevolumeyear)\hfill\fi}
\def\@evenfoot{\@oddfoot}
\makeatother
\fi

\newwrite\gtoutfile
\long\gdef\makeheadfile{  
{\def\\{, }\def\s{ }
\immediate\openout\gtoutfile head.xxx
\immediate\write\gtoutfile{To: math@arxiv.org}
\immediate\write\gtoutfile{Subject: put or rep NNNNN:pppp}
\immediate\write\gtoutfile{--text follows this line--}
\immediate\write\gtoutfile{Proxy-for: \ifx\theasciiauthors\relax
\theauthors\else\theasciiauthors\fi\s<\ifx\theasciiemail\relax\theemail\else\theasciiemail\fi>}
\immediate\write\gtoutfile{\noexpand\\}
\immediate\write\gtoutfile{Authors: \ifx\theasciiauthors\relax
\theauthors\else\theasciiauthors\fi}
\immediate\write\gtoutfile{Title: \ifx\theasciititle\relax
\thetitle\else\theasciititle\fi}
\immediate\write\gtoutfile{Subj-class: GT or SG or MG etc}
\immediate\write\gtoutfile{MSC-class: \theprimaryclass\ifx\thesecondaryclass\relax\else, \thesecondaryclass\fi}
\immediate\write\gtoutfile{Journal-ref: Geom. Topol. \thevolumenumber
(\thevolumeyear) \startpage-\finishpage}
\immediate\write\gtoutfile{Comments: Published by Geometry and Topology at}
\immediate\write\gtoutfile{\s\s http://www.maths.warwick.ac.uk/gt/GTVol\thevolumenumber/paper\thepapernumber.abs.html}
\immediate\write\gtoutfile{\noexpand\\}
\immediate\write\gtoutfile{}
\ifx\theasciiabstract\relax
\immediate\write\gtoutfile{\theabstract}\else
\immediate\write\gtoutfile{\theasciiabstract}\fi
\immediate\write\gtoutfile{}
\immediate\write\gtoutfile{\noexpand\\}
\immediate\write\gtoutfile{}
\immediate\closeout\gtoutfile}}  

\def\maketitlepage{\maketitlep\makeheadfile}
\let\maketitle\maketitlepage

\def\ifplaintex{\expandafter\ifx\csname documentclass\endcsname\relax}

\ifplaintex 
\else
\fi

\expandafter\ifx\csname beginpicture\endcsname\relax
\expandafter\ifx\csname documentclass\endcsname\relax
\input pictex \else
\input prepictex \input pictex \input postpictex \fi\fi

\def\gt{{\mathsurround=0pt\it $\cal G\mskip-2mu$eometry \&\ 
$\cal T\!\!$opology}}        

\def\gtp{{\mathsurround=0pt\it $\cal G\mskip-2mu$eometry \&\ 
$\cal T\!\!$opology $\cal P\!$ublications}}  

\def\lognumber#1{\def\thelognumber{#1}}
\def\volumenumber#1{\def\thevolumenumber{#1}}
\def\papernumber#1{\def\thepapernumber{#1}}
\def\volumeyear#1{\def\thevolumeyear{#1}}

\def\pagenumbers#1#2{\def\startpage{#1}\def\finishpage{#2}}
\def\published#1{\def\publishdate{#1}}
\def\proposed#1{\def\theproposer{#1}}
\def\seconded#1{\def\theseconders{#1}}
\def\received#1{\def\receiveddate{#1}}

\def\accepted#1{\def\accepteddate{#1}}

\def\asciiaddress#1{\def\theasciiaddress{#1}}

\long\def\asciiabstract#1{\long\def\theasciiabstract{#1}}
\def\asciikeywords#1{\def\theasciikeywords{#1}}

\let\\\par\let\thelognumber\relax
\let\thevolumenumber\relax\let\thepapernumber\relax
\let\thevolumeyear\relax\let\thesamplenumber\relax\let\startpage\relax
\let\finishpage\relax\let\publishdate\relax\let\receiveddate\relax
\let\reviseddate\relax\let\accepteddate\relax\let\theasciititle\relax
\let\theasciiauthors\relax\let\theasciiaddress\relax
\let\theasciiabstract\relax\let\theasciikeywords\relax
\let\theasciiemail\relax\let\theshortauthors\relax\let\theshorttitle\relax

\long\def\maketitlep{   

\count0=\startpage

\gt\hfill      
\beginpicture
\setcoordinatesystem units <0.33truein, 0.33truein> point at 2.2 0.9
\setplotsymbol ({$\cal G$})
\plotsymbolspacing=9truept
\circulararc 315 degrees from 0 1 center at 0 0
\setplotsymbol ({$\cal T$})
\circulararc 315 degrees from 1 -1 center at 1 0
\endpicture
\break
{\small\ifx\thesamplenumber\relax 
Volume \else Sample
\fi\thevolumenumber\ (\thevolumeyear)
\startpage--\finishpage\nl
Published: \publishdate}
\vglue 0.5truein plus 0.4fil minus 0.1truein

{\parskip=0pt\leftskip 0pt plus 1fil\def\\{\par\smallskip}{\ifplaintex\large
\else\Large\fi\bf\thetitle}\par\medskip}   

\vglue 0pt plus 0.1fil 

{\parskip=0pt\leftskip 0pt plus 1fil\def\\{\par}{\sc\theauthors}
\par\medskip}

\vglue 0pt plus 0.1fil 

{\small\parskip=0pt\let\newline\\
{\leftskip 0pt plus 1fil\def\\{\par}{\sl\theaddress}\par}
\expandafter\ifx\theemail\relax    
\relax\else\vglue 5pt plus 0.02fil minus 2pt\def\\{\stdspace{\rm 
and}\stdspace} 
\cl{Email:\stdspace\tt\theemail}\fi
\ifx\theurl\relax                  
\relax\else\vglue 5pt plus 0.02fil minus 2pt\def\\{\stdspace{\rm 
and}\stdspace}
\cl{URL:\stdspace\tt\theurl}\fi\par}

\vglue 7pt plus 0.3fil minus 3pt

{\bf Abstract}
\vglue 5pt plus 0.1fil minus 2pt

\theabstract

\vglue 7pt plus 0.3fil minus 3pt

{\bf AMS Classification numbers}\quad Primary:\quad \theprimaryclass

Secondary:\quad \thesecondaryclass

\vglue 5pt plus 0.3fil minus 2pt

{\bf Keywords}\quad \thekeywords

\vglue 10pt plus 0.5fil minus 5pt

{\small  Proposed: \theproposer\hfill Received: \receiveddate\nl
Seconded: \theseconders\hfill 
\ifx\reviseddate\relax                         
Accepted: \accepteddate                        
\else
Revised: \reviseddate                          
\fi}
\eject
}       

\let\maketitlepage\maketitlep
\let\maketitle\maketitlepage

\font\phead=cmsl9 scaled 950
\font=cmsl9 scaled 1050
\font\pnum=cmbx10 scaled 913
\font\lnum=cmbx10 
\font\pfoot=cmsl9 scaled 950
\font=cmsl9 scaled 1050
\ifplaintex
\headline{\vbox to 0pt{\vskip -4.5mm\line{\small\phead\ifnum
\count0=\startpage ISSN 1364-0380 (on line)
1465-3060 (printed) \hfill {\pnum\folio}\else\ifodd\count0\def\\{ }%
\ifx\theshorttitle\relax\thetitle\else\theshorttitle\fi\hfill{\pnum\folio}
\else\def\\{ and }{\pnum\folio}\hfill\ifx\theshortauthors\relax\theauthors
\else\theshortauthors\fi\fi\fi}\vss}}
\footline{\vbox to 0pt{\vglue 0mm\line{\small\pfoot\ifnum\count0=\startpage
\copyright\ \gtp\hfill\else
\gt, Volume \thevolumenumber\ (\thevolumeyear)\hfill\fi}\vss
}}
\else
\makeatletter
\def\@oddhead{{\small\lhead\ifnum\count0=\startpage ISSN 1364-0380 (on line)
1465-3060 (printed) \hfill {\lnum\number\count0}\else\ifodd\count0
\def\\{ }\ifx\theshorttitle\relax \thetitle \else\theshorttitle\fi\hfill
{\lnum\number\count0}\else\def\\{ and }{\lnum\number\count0}
\hfill\ifx\theshortauthors\relax 
\theauthors\else\theshortauthors\fi\fi\fi}}\def\@evenhead{\@oddhead}
\def\@oddfoot{\small\lfoot\ifnum\count0=\startpage\copyright\ \gtp\hfill\else
\gt, Volume \thevolumenumber\ (\thevolumeyear)\hfill\fi}
\def\@evenfoot{\@oddfoot}
\makeatother
\fi

\newwrite\gtoutfile
\long\gdef\makeheadfile{  
{\def\\{, }\def\s{ }
\immediate\openout\gtoutfile head.xxx
\immediate\write\gtoutfile{To: math@arxiv.org}
\immediate\write\gtoutfile{Subject: put or rep NNNNN:pppp}
\immediate\write\gtoutfile{--text follows this line--}
\immediate\write\gtoutfile{Proxy-for: \ifx\theasciiauthors\relax
\theauthors\else\theasciiauthors\fi\s<\ifx\theasciiemail\relax\theemail\else\theasciiemail\fi>}
\immediate\write\gtoutfile{\noexpand\\}
\immediate\write\gtoutfile{Authors: \ifx\theasciiauthors\relax
\theauthors\else\theasciiauthors\fi}
\immediate\write\gtoutfile{Title: \ifx\theasciititle\relax
\thetitle\else\theasciititle\fi}
\immediate\write\gtoutfile{Subj-class: GT or SG or MG etc}
\immediate\write\gtoutfile{MSC-class: \theprimaryclass\ifx\thesecondaryclass\relax\else, \thesecondaryclass\fi}
\immediate\write\gtoutfile{Journal-ref: Geom. Topol. \thevolumenumber
(\thevolumeyear) \startpage-\finishpage}
\immediate\write\gtoutfile{Comments: Published by Geometry and Topology at}
\immediate\write\gtoutfile{\s\s http://www.maths.warwick.ac.uk/gt/GTVol\thevolumenumber/paper\thepapernumber.abs.html}
\immediate\write\gtoutfile{\noexpand\\}
\immediate\write\gtoutfile{}
\ifx\theasciiabstract\relax
\immediate\write\gtoutfile{\theabstract}\else
\immediate\write\gtoutfile{\theasciiabstract}\fi
\immediate\write\gtoutfile{}
\immediate\write\gtoutfile{\noexpand\\}
\immediate\write\gtoutfile{}
\immediate\closeout\gtoutfile}}  

\def\maketitlepage{\maketitlep\makeheadfile}
\let\maketitle\maketitlepage